\documentclass[12pt]{amsart}

\usepackage{amsfonts,amssymb, latexsym, epsfig, epic,enumitem}
\usepackage{amsmath,amsthm,lipsum,commath,graphicx,setspace}
\usepackage{tikz}
\usepackage{verbatim,color,xcolor}
\usetikzlibrary{arrows,positioning}
\graphicspath{{figures/}}

\newtheorem{defn}{Definition}[section]
\newtheorem{thm}{Theorem}[section]
\theoremstyle{plain}

\newtheorem{lemma}{Lemma}[section]
\newtheorem{cor}{Corollary}[section]
\newtheorem{prop}{Proposition}[section]

\begin{document}

\title{Belted sum decompositions of fully augmented links}

\author[P. Morgan]{Porter Morgan}

\author[B. Ransom]{Brian Ransom}

\author[D. Spyropoulos]{Dean Spyropoulos}

\author[C. Ziegler]{Cameron Ziegler}

\author[R. Trapp]{Rolland Trapp}
\address{California State University, San Bernardino\\Dept. of Mathematics\\
5500 University Pkwy\\
San Bernardino, CA 92407}
\email{rtrapp@csusb.edu}

\begin{abstract}
Given two orientable, cusped hyperbolic 3-manifolds containing certain thrice-punctured spheres, Adams gave a diagrammatic definition for a third such manifold, their belted sum. Fully augmented links, or FALs, are hyperbolic links constructed by augmenting a link diagram. This work considers belted sum decompositions in which all manifolds involved are FAL complements.  To do so, we provide explicit classifications of thrice-punctured spheres in FAL complements, making them easily recognizable.  These classifications are used to characterize belted sum prime FALs geometrically, combinatorially and diagrammatically. Finally we prove that, in the context of belted sums, every FAL complement canonically decomposes into FALs which are either prime or two-fold covers of the Whitehead link.
\end{abstract}

\maketitle
\section{Introduction}\label{sec:Intro}

\begin{figure}[h]

\begin{center}
\includegraphics{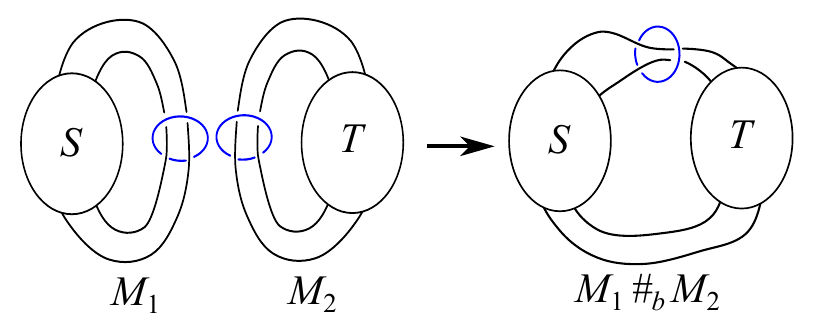}
\end{center}

\caption{Forming the belted sum of $M_1$ and $M_2$}
\label{fig:BS}
\end{figure}

Within the context of orientable hyperbolic three-manifolds, Colin Adams defined a natural belted-sum operation (see \cite{ad1}).  Adams uses embedded thrice-punctured spheres to form a sum of two cusped, orientable, hyperbolic three-manifolds as follows.  First, slice along thrice-punctured spheres in each, resulting in manifolds whose boundaries contain two thrice-punctured spheres. Then glue the manifolds together along the boundary thrice-punctured spheres.  Adams proves that the result is again hyperbolic and that volumes add under this operation.  A belted sum is a particular type of this operation, which Adams defines diagrammatically as in Figure \ref{fig:BS}.

This paper studies the reverse of the belted sum operation on a particular class of hyperbolic links called fully augmented links, or FALs.  Hyperbolic FALs are constructed from certain link diagrams by augmenting each twist region with an unknot--placing a trivial knot around each twist region and removing all full twists (see Figure \ref{fig:FA}, and Section \ref{sec:FAL} for more detail).  FALs were introduced by Agol and Thurston in the appendix of \cite{la}, where they used FALs to improve Lackenby's upper bound on volumes of hyperbolic alternating link complements.  Since their introduction, FALs have proven useful in a variety of other contexts as well. The geometry of an FAL complement can be described by gluing a right-angled ideal polyhedron to its reflection using products of reflections in its faces (see Section \ref{sec:FAL} or Purcell's paper \cite{pu1}).  In addition, most Dehn fillings of the augmenting unknot components of an FAL yield hyperbolic links whose geometry can be related to that of the FAL. For example, FALs have been used to produce links that admit no exceptional surgeries, to bound the genera of surfaces in highly twisted link complements, to create diagrammatic bounds on volumes of sufficiently twisted links, and to bound crosscap numbers of alternating links (see \cite{fp}, \cite{bft}, \cite{fkp}, and \cite{kl}, respectively).  Thus FALs are a convenient and useful class of hyperbolic links to study.

In this work we focus on decomposing a given FAL complement into the belted sum of two simpler FALs, where term ``simpler" will develop to have meaning in topological, combinatorial and geometric contexts.  An FAL which cannot be so decomposed is called $b$-prime and, since several properties of manifolds are well-behaved under belted-sums, a better understanding of $b$-prime FALs will shed light on FALs in general.  We now highlight two potential applications.

First, since volumes are additive under belted sums, understanding volumes of $b$-prime FALs can be applied to the study of FAL volumes in general.  Volumes of FALs, in turn, are of interest since they are related to (conjectured) minimal volume $n$-cusped hyperbolic manifolds. Yoshida has shown that an FAL is the minimal volume orientable 4-cusped hyperbolic manifold (the fully twisted Borromean rings of Figure \ref{fig:Tet8}). Moreover, Agol \cite{A2010} has conjectured that coverings of the Whitehead link are minimal volume for links with at least $11$ components.  Purcell proves (see \cite[Proposition 3.6]{pu1}) that $2(c-1)v_8$ is a sharp lower bound for the volume of an FAL with $c$ crossing circles, where $v_8$ is the volume of a regular ideal octahedron, and there are octahedral FALs that satisfy Agol's conjectured minimal volume.  These observations indicate that knowledge of FAL volumes, combined with the tractable geometry of FALs, could result in further applications. 

Arithmetic invariants and, in particular, invariant trace fields, behave nicely under belted sum operations, providing a second application for $b$-prime FALs.  Indeed, Section 5.6 of Maclachlan and Reid (see \cite[Theorem 5.6.1]{mr}) shows that the invariant trace field of $M_1 \#_b M_2$ is the compositum of the fields of the summands $M_i$.  Moreover, several authors have used the tractable geometry of FALs to study questions of arithmeticity and commensurability. Meyer-Millichap-Trapp computed invariant trace fields of minimally twisted chain links with an even number of components, and determined which such links are arithmetic and which are commensurable (see \cite{MeMiTr2020}). In addition, Rochy Flint (\cite{rf}) developed techniques for calculating invariant trace fields of FALs.  Thus the explicit geometry of $b$-prime FALs combined with the belted sum operation can provide insight into invariant trace fields, arithmeticity and commensurability questions for FALs.

It behooves us, then, to analyze belted sum decompositions of FALs. The first step in this analysis is to classify embedded, totally geodesic thrice-punctured spheres in FAL complements.  Every FAL complement admits an orientation-reversing involution whose fixed point set is an embedded, totally geodesic \emph{reflection surface}. We will see, in Lemma \ref{lem:Standard}, that thrice-punctured sphere components of a reflection surface are never part of a belted sum decomposition.  For this reason we focus on thrice-punctured spheres not contained in the reflection surface.  Section \ref{sec:3ps} develops the necessary tools to classify such \emph{non-reflection} thrice-punctured spheres, showing there are three types.  The classification is given in Theorem \ref{thm:NR3ps}, stated below. The three types of spheres are illustrated in Figures \ref{fig:NDisks} and \ref{fig:SingSepDisk} (see also Definitions \ref{defn:NDisks} and \ref{defn:SDisk}).

\vspace{6pt}
\noindent
\textbf{Theorem \ref{thm:NR3ps}.} \emph{Let $S$ be a non-reflection thrice-punctured sphere in an FAL complement $M = \mathbb{S}^3\setminus \mathcal{A}$.  If $\mathcal{A}$ is not the fully twisted Borromean rings, then $S$ is orthogonal to the reflection surface in $M$ and is either a crossing, longitudinal or singly-separated disk.}
\vspace{6pt}

The three types--crossing disks, longitudinal disks, singly-separated disks--of non-reflection thrice-punctured spheres are classified by their punctures and how they intersect the reflection surface.  Thus their geometry, and how it relates to the overall geometry of the FAL complement, is made quite explicit (see Propositions \ref{prop:3psDisjointFromCD} and \ref{prop:3psSeparatesD}).  Moreover, they are easy to recognize using the tools developed in Section \ref{sec:FAL}.

Section \ref{sec:SP} determines which thrice-punctured spheres of Theorem \ref{thm:NR3ps} can be used in a belted sum decomposition.  Section \ref{sec:FALPrime} applies results of Section \ref{sec:SP} to give several characterizations of $b$-prime FALs.  The background information of Section \ref{sec:FAL} is required to describe two of the characterizations given, so we focus on the third.  There is a procedure to construct an FAL from a link diagram and FAL complements are (essentially) $b$-prime if and only if the diagram has only trivial flype orbits.  The result is Theorem \ref{thm:PrimeIsTrivialOrbit}, stated below (technical definitions are postponed until Section \ref{sec:FALPrime}).  

\vspace{6pt}
\noindent
\textbf{Theorem \ref{thm:PrimeIsTrivialOrbit}.} \emph{Let $D$ be a twist-reduced, non-split, prime link diagram with more than two twist regions, and let $\mathcal{A}$ be a full augmentation of $D$.  
The fully augmented link $\mathcal{A}$ is $b$-prime if and only if it induces only trivial flype orbits.}
\vspace{6pt}

We briefly mention, without detail, two other characterizations of $b$-primality, and refer the reader to Section \ref{sec:FALPrime} for precise statements. Theorem \ref{thm:CrushCharPrime} is a combinatorial analogue of Theorem \ref{thm:PrimeIsTrivialOrbit}, and is stated in terms of perfect matchings on cubic planar graphs associated with FALs.  Corollary \ref{cor:PrimeIs1Disk} provides a geometric characterization of $b$-prime FALs stated in terms of the absence of certain thrice-punctured spheres. The combinatorial characterization provides an efficient method for tabulating $b$-prime FALs by number of crossing circles. All three characterizations of $b$-primality provide concrete approaches to the geometric and arithmetic questions described above.

An immediate application of the geometric characterization of Corollary \ref{cor:PrimeIs1Disk} is the definition of a canonical belted sum decomposition in Section \ref{sec:CanonDec}.  Of the three types of non-reflection thrice-punctured spheres (crossing, longitudinal, and singly-separated disks), the canonical decomposition is obtained by decomposing only along crossing disks (see Definition \ref{defn:CanonicalBS}).  The result, stated here in Theorem \ref{thm:CanonIsBPrime}, is that the summands of the canonical decomposition are either $b$-prime or 2-fold covers of the Whitehead link.

\vspace{6pt}
\noindent
\textbf{Theorem \ref{thm:CanonIsBPrime}.} \emph{Let $\mathcal{A}$ be an FAL with at least three crossing circles and let $M=\mathbb{S}^3\setminus\mathcal{A}$ be its complement. Each summand in the canonical belted sum decomposition of $M$ is either a $b$-prime FAL complement or the complement of the Borromean rings with at least one flat crossing circle.
}
\vspace{6pt}

The paper is organized as follows.  Section \ref{sec:FAL} reviews the definition of FALs, and the tools associated with them necessary for our study.  Section \ref{sec:3ps} characterizes thrice-punctured spheres in FAL complements, and Section \ref{sec:SP} classifies which can be used in belted-sum decompositions. Section \ref{sec:FALPrime} gives the three characterizations of $b$-primality in FAL complements, and Section \ref{sec:CanonDec} introduces the canonical belted sum decomposition of an FAL.

\section{Fully Augmented Links} \label{sec:FAL}

\begin{figure}[h]
\[
\begin{array}{ccc}
\includegraphics[width=1.4in]{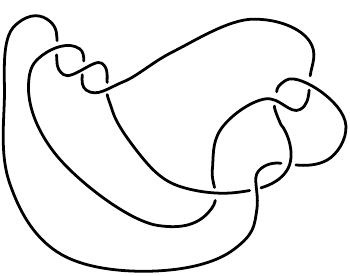}&\includegraphics[width=1.4in]{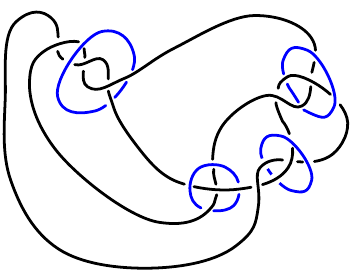}& \includegraphics[width=1.4in]{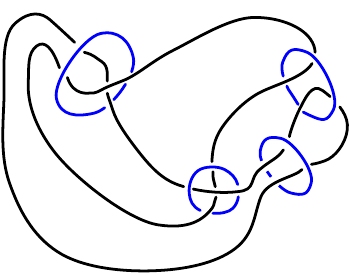}\\
(a) \textrm{ Link diagram} &(b) \textrm{ Augmented }& (c) \textrm{ Fully augmented}
\end{array}
\]
\caption{Fully augmenting a link}
\label{fig:FA}
\end{figure}

Augmenting a link consists of placing a trivial component around a twist of two strands.  Adams, in \cite{ad2}, showed that augmenting an alternating hyperbolic link results in another hyperbolic link.  Fully augmented links (or FALs) are the result of augmenting every twist region and removing all full twists (see Figure \ref{fig:FA}).  Purcell relaxed the restriction that a diagram be alternating to show that an FAL created from any prime, twist reduced diagram with at least two twist regions is hyperbolic (see \cite[Theorem 6.1]{pu3}).

This section reviews both geometric and combinatorial approaches toward studying FAL complements.  First, slicing an FAL complement along the two-cells of a standard cell decomposition facilitates the transition from topology to geometry, yielding a cell decomposition of the two-sphere.  Andreev's Theorem implies that this cell decomposition corresponds to a circle packing of $\mathbb{S}^2$ (see, for example, \cite[Corollary 13.6.2]{th}).  Geometrically the circle packing is the ``footprint", if you will, of a right-angled ideal polyhedron in hyperbolic space which is unique up to M\"obius transformation. The original FAL complement, then, is made up of two such polyhedra.  There are also two graphs naturally associated with the circle packing--the nerve and its dual--which provide convenient combinatorial descriptions of FALs.  Geometric and combinatorial tools will be used extensively throughout this paper, and we proceed with a more thorough description of them.

Let $\mathcal{A}_L$ be the fully augmented link arising from $L$ as described above.  The additional trivial components are \emph{crossing circles}, and the components that remain from the original link $L$ will be called \emph{knot circles}.  A crossing circle with a half-twist is called \emph{twisted}, while those without are \emph{flat}.  Each crossing circle bounds a disk, called a \emph{crossing disk}, which is punctured twice by knot circles.  Crossing circles can bound more than one crossing disk, and when this occurs it always results in a possible belted-sum decomposition (see Theorem \ref{thm:GeneralBS}).  In what follows we assume a specific crossing disk is chosen for each crossing circle.  Moreover, the complement $M = S^3\setminus \mathcal{A}_L$ admits reflective symmetry across the projection plane, followed by a Dehn twist on twisted crossing circles.  Purcell shows that the fixed point set of this reflection is an embedded, totally geodesic surface $R\subset M$, called the \emph{reflection surface} (see \cite[Lemma 2.1]{pu1}). If there is more than one reflection surface in an FAL complement we assume a specific choice of reflection surface has been made, just as a crossing disk is chosen for each crossing circle.

\begin{figure}[h]
\[
\begin{array}{ccc}
\includegraphics[width=1.4in]{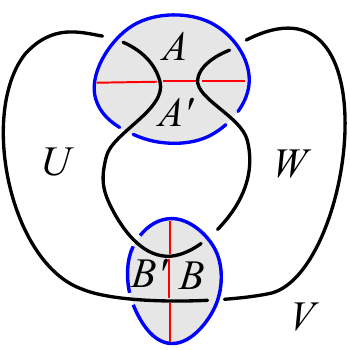}&\includegraphics[width=1.4in]{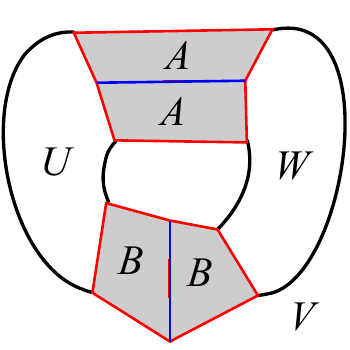}&\includegraphics[width=1.4in]{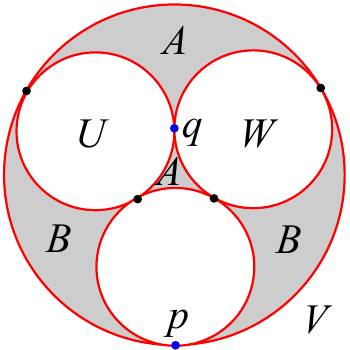}\\
(a) \textrm{ Standard cell } & (b) \textrm{ Sliced 2-cells } & (c)\textrm{ Circle packing}\\
\textrm{  decomposition}& \textrm{ in }\partial B_+ & \textrm{ for }P_+
\end{array}
\]
\caption{From cell decomposition to circle packing}
\label{fig:CellDecomp}
\end{figure}

The complement $M = S^3\setminus \mathcal{A}_L$ of every fully augmented link $\mathcal{A}_L$ admits a \emph{standard cell decomposition} $\mathcal{C}$.  The standard cell decomposition is most easily seen when all crossing circles are flat, so we begin with this case.  There are two kinds of 2-cells in $\mathcal{C}$: the regions of the projection plane (e.g. $U$ and $V$ in Figure \ref{fig:CellDecomp}$(a)$) and the crossing disks (so $\mathcal{C}$ depends on our initial choice of crossing disks).  We refer to the projection plane 2-cells as \emph{reflection} 2-cells since reflection across the plane of projection is a symmetry of the link $\mathcal{A}_L$.  The 2-cells from crossing disks are \emph{crossing} 2-cells.  The projection plane cuts through the middle of each crossing disk, and the curves of intersection are the 1-cells of the standard cell decomposition (see the vertical and horizontal line segments of Figure \ref{fig:CellDecomp}$(a)$).  Crossing disks are cut into two triangles by the 1-cells, one above and one below the plane of projection (e.g. $A, A'$ in Figure \ref{fig:CellDecomp}$(a)$) .

There are no 0-cells in the standard cell decomposition as all 1-cells have endpoints on $\mathcal{A}_L$, which will become ideal points. The two 3-cells $B_{\pm}^3$ correspond to the regions of $M$ above and below the projection plane.  

Slicing along the reflection 2-cells separates $M$ into two pieces which are reflections of each other, and crossing disks contribute a triangle in each half.  Further slicing along the crossing disks yields the two three balls $B_{\pm}^3$, and the boundary of $B_+$ (viewed from ``inside" $B_+$) is pictured in Figure \ref{fig:CellDecomp}$(b)$.  Shrinking the arcs of the original link to vertices gives a cell decomposition on the boundary of $B_+$. This cell decomposition can be checkerboard colored so that the unshaded faces correspond to reflection 2-cells, while shaded faces are all triangular and come from crossing disks.  

The unshaded faces can be isotoped to form a circle packing by Andreev's Theorem (see \cite{an}), where two circles are tangent in the circle packing if the corresponding regions in the plane share a knot circle or are both punctured by the same crossing circle. Thus the circles corresponding to regions $U$ and $W$ are tangent at the vertex $q$ corresponding to the crossing circle bounding $A$ (see Figure \ref{fig:CellDecomp}$(c)$).  

This circle packing leads directly to a geometric description of an FAL complement.  A further consequence of Andreev's Theorem is that $P_{\pm}$ can be realized as right-angled ideal polyhedra which are unique up to M\"obius transformation (see \cite{pu1} and \cite[Corollary 13.6.2]{th} for more detail).  We refer to the polyhedra $P_{\pm}$ as \emph{standard polyhedra} for the FAL $\mathcal{A}_L$.

\begin{figure}[h]
\[
\begin{array}{ccc}
\includegraphics[width=1.2in]{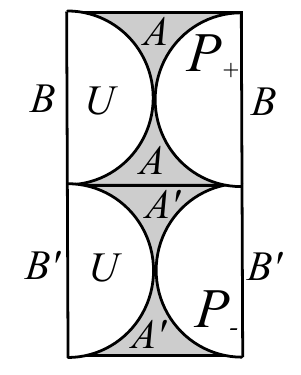}&\includegraphics[width=1.2in]{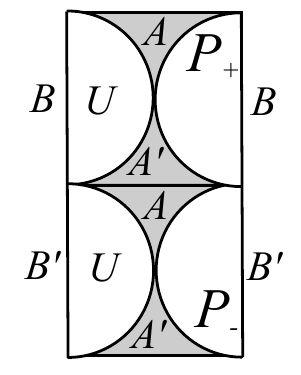}& \includegraphics[width=1.2in]{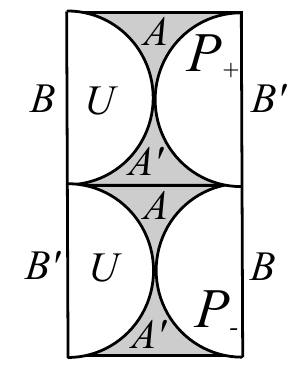}\\
(a) \textrm{ Borromean Rings } &(b) \textrm{ Once twisted }& (c) \textrm{ Fully twisted}\\
&\phantom{(b)}\textrm{ Borromean Rings }&\phantom{(c)}\textrm{ Borromean Rings }
\end{array}
\]
\caption{Comparing flat and twisted gluing patterns}
\label{fig:F}
\end{figure}

Since the 3-cells $B_{\pm}^3$ are reflections of each other across the projection plane, the circle packings on their boundaries lead to polyhedra $P_{\pm}$ that are reflections of each other.  To find a fundamental domain for $M$, let $P_-$ be the reflection of $P_+$ across an unshaded face $G$.  Then $\mathcal{F} = P_+\cup P_-$ is a fundamental domain for $M$.  Figure \ref{fig:F}$(a)$ illustrates the fundamental domain $\mathcal{F}$ for the Borromean rings of Figure \ref{fig:CellDecomp}$(a)$ after transforming the vertex $p$ to infinity.

The gluing maps on the faces of $\mathcal{F}$ are quite explicit.  To describe the gluings between unshaded faces we refer to faces, edges and vertices of $\mathcal{F}$ that are reflections of each other as \emph{corresponding} faces, edges or vertices, respectively. If $F$ and $F'$ are corresponding unshaded faces, then they are glued by identifying corresponding edges and vertices.  This is realized by the isometry $r_G\circ r_F : F\to F'$, where $r_F$ and $r_G$ represent reflections across the respective faces (recall that $G$ is the shared face of $P_{\pm}$).  Indeed, reflection in $F$ fixes $F$ point-wise while reflection in $G$ identifies corresponding components of $F$ and $F'$.

For a flat FAL, each shaded triangle on $P_+$ will be glued to an adjacent shaded triangle by the parabolic isometry that fixes their shared point (for example, the point $q$ between the triangles labeled $A$ in Figure \ref{fig:CellDecomp}).  Similarly for shaded triangles on $P_-$.  

The universal cover $\widetilde{M}$ of $M$ is tessellated by the images of $\mathcal{F}$ under the group generated by the gluing maps just described.  The fundamental domain $\mathcal{F}$ is an example of what we'll call a standard domain for $M$, since it was formed using its standard cell decomposition.  A \emph{standard domain} for $M$ is the union of $P_+$ with any copy of $P_-$ sharing an unshaded face.  

To see how twisted crossing circles change the above description, note that a twisted crossing circle can be obtained from a flat one by slicing along the crossing disk and regluing with a half-twist.  This glues the top and bottom of one crossing disk to the bottom and top of the other.  Carrying the projection plane along in the process yields the reflection 2-cells which form a reflection surface for the twisted disk.  The reflection is the homeomorphism defined by reflecting in the plane of projection followed by full twists inside each twisted crossing circle.  

FALs that differ by half-twists will be called \emph{half-twist partners}.  Note that the polyhedra $P_{\pm}$ are the same for the half-twist partners, but shaded triangles in opposite polyhedra are identified.  Figure \ref{fig:F}$(a)$, for example, illustrates a fundamental domain for the Borromean rings in which the vertex $p$ of Figure \ref{fig:CellDecomp}$(c)$ is placed at infinity in the upper half-space model.  Figure \ref{fig:F} parts $(b)$ and $(c)$ illustrate the gluing patterns for the Borromean rings' half-twist partners.

\begin{figure}[h]
\[
\begin{array}{ccc}
\includegraphics[width=1.4in]{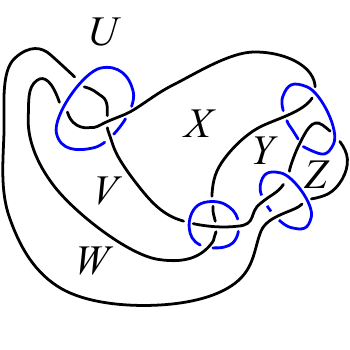}&\includegraphics[width=1.4in]{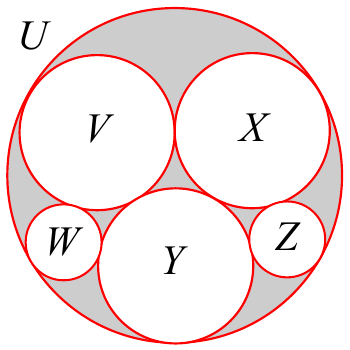}& \includegraphics[width=1.4in]{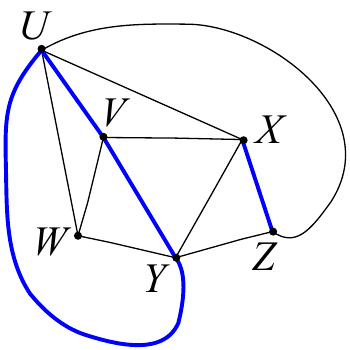}\\
(a) \textrm{ An FAL } &(b) \textrm{ Circle packing }& (c) \textrm{ Colored nerve }\Gamma
\end{array}
\]
\caption{Colored nerve of an FAL}
\label{fig:Nerve}
\end{figure}

The standard polyhedra $P_{\pm}$ give a nice geometric description of $M$ and, following \cite{pu1}, we describe combinatorial descriptions that also arise from the circle packing.  The \emph{nerve} of a circle packing is the planar graph obtained by placing a vertex at the center of each circle and an edge through points of tangency. Since the shaded regions are triangular, the nerve of an FAL circle packing is always a triangulation of $\mathbb{S}^2$.  Figure \ref{fig:Nerve} illustrates a FAL with reflection 2-cells labeled, the circle packing for $P_+$, and its corresponding nerve. 

The nerve $\Gamma$ can be painted, as described in \cite{pu1}, to record which shaded triangles are glued or, equivalently, which ideal vertices correspond to crossing circle cusps. Glued shaded triangles must share an ideal vertex that corresponds to a crossing circle cusp, so simply paint edges in the nerve through these shared vertices as in Figure \ref{fig:Nerve}$(c)$.  Painting the nerve in this way results in a triangulation of $\mathbb{S}^2$ in which one edge of each triangle is painted.  

Conversely, given a triangulation of $\mathbb{S}^2$ for which each triangle has one painted edge, one can construct a unique flat FAL from it.  Each painted edge is replaced with a crossing circle that punctures the plane near its vertices.  Knot circles lie in the plane and connect midpoints of non-painted edges through the crossing circles just constructed.  The result is a flat FAL.  Of course, arbitrarily twisting the crossing circles leads to other FALs that give rise to the same painted nerve.

\begin{figure}[h]
\begin{center}
\includegraphics[width=1.4in]{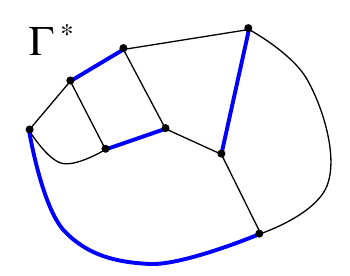}
\end{center}
\caption{Painted crushtacean $\Gamma^{\ast}$}
\label{fig:Crush}
\end{figure}

The planar dual $\Gamma^{\ast}$ of the nerve $\Gamma$ provides a second combinatorial structure for studying FALs. Following Chesebro, DeBlois and Wilton, we call the planar dual $\Gamma^{\ast}$ the \emph{crushtacean} (see \cite{CDW2012} and Figure \ref{fig:Crush}). Vertices of $\Gamma^{\ast}$ correspond to shaded triangles in the circle packing and edges to points of tangency between them.

Since the nerve $\Gamma$ is a triangulation of $\mathbb{S}^2$, the crushtacean $\Gamma^{\ast}$ is a trivalent planar graph.  A painting on $\Gamma$ determines one on the crushtacean simply by painting edges that are dual to painted edges of $\Gamma$.  Note that each vertex of $\Gamma^{\ast}$ has one painted edge since the same is true of each face of $\Gamma$, so a painted crushtacean is simply a perfect matching on a trivalent planar graph.  Moreover, vertices sharing a painted edge correspond to shaded triangles in $P_+$ that are glued. 

\begin{figure}[h]
\begin{center}
\includegraphics{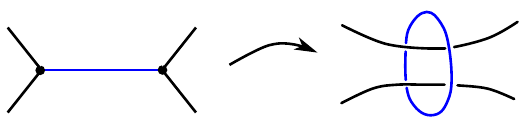}
\end{center}
\caption{Creating a FAL from a painted crushtacean}
\label{fig:FALFromCrush}
\end{figure}

A painted nerve and crushtacean, then, correspond to each FAL.  Conversely Purcell shows in \cite[Lemma 2.4]{pu1} how to create an FAL from a painted nerve $\Gamma$.  Although \cite[Lemma 2.4]{pu1} is stated in terms of painted nerves, the proof's construction resorts to crushtaceans.  Given the painted crushtacean $\Gamma^{\ast}$, simply replace each painted edge with a crossing circle that links it once and smooth the adjacent unpainted edges as in Figure \ref{fig:FALFromCrush} (which is Figure 5 of \cite{pu1}).  The FAL resulting from this process determines the same polyhedra with appropriate gluing pattern. Of course, one alters this process appropriately to accommodate half-twists partners with the same painted crushtacean $\Gamma^{\ast}$.

Features of the circle packing and crushtacean will be useful in characterizing thrice-punctured spheres and belted sums of FAL complements.  We  introduce two that will be particularly useful.  

A shaded triangular interstice in the circle packing corresponds to a shaded triangular face of a standard polyhedron, say $P_+$.  The shaded face is an ideal triangle which orthogonal to the unshaded faces of $P_+$ that it intersects, and which projects to half a crossing disk in the FAL complement.  Another way to describe this is that three mutually tangent circles that bound an interstice correspond to an ideal triangular face of $P_+$ that projects to (a portion of) a crossing disk.  More generally, we now observe that three mutually tangent circles in the circle packing correspond to an ideal triangle in $P_{\pm}$ that intersects unshaded faces orthogonally.  Later, Proposition \ref{prop:3psDisjointFromCD} will demonstrate that these always project to thrice-punctured spheres in the FAL complement.

\begin{figure}[h]
\[
\begin{array}{ccc}
\includegraphics{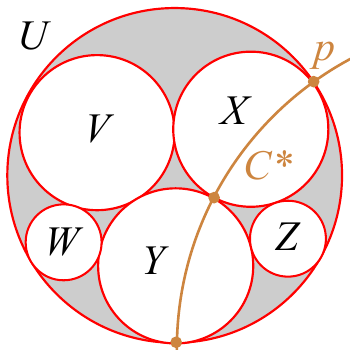}&\includegraphics{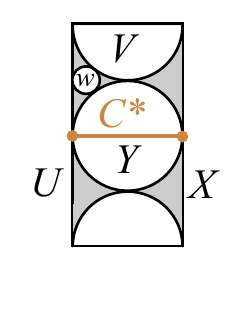}& \includegraphics{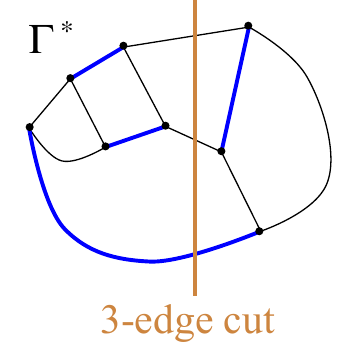}\\
(a) \textrm{ Non-standard triple}&(b) \textrm{ Vertex $p$ at infinity}& (c) \textrm{ Non-trivial cut}
\end{array}
\]
\caption{A non-standard triple and corresponding non-trivial 3-edge cut}
\label{fig:NSMT}
\end{figure}

Let $C_1, C_2, C_3$ be three mutually tangent circles in the circle packing, let $p,q,r$ denote their points of tangency, and let $C^{\ast}$ the circle determined by $p,q,r$.  We first show that the hyperbolic plane bounded by $C^{\ast}$ is orthogonal to those determined by the $C_i$. Apply a M\"obius transformation to the four circles that takes $p$ to infinity and use the same labels to refer to the transformed figures.  Then $C_1$ and $C_2$ are parallel lines and $C_3$ is a circle tangent to each at $q,r$, so that segment $qr$ is a diameter of $C_3$.  The transformed $C^{\ast}$, then, is the line $\overleftrightarrow{qr}$, which is orthogonal to the $C_i$ and the original $C^{\ast}$ is orthogonal to the original $C_i$.  Thus $C^{\ast}$ bounds a plane $\mathcal{P}\subset\mathbb{H}^3$ and $\mathcal{P}\cap P_{\pm}$ is an ideal triangle orthogonal to unshaded faces of $P_{\pm}$.  If the $C_i$ bound a shaded triangular interstice in the circle packing, then $\mathcal{P}\cap P_{\pm}$ is a shaded face of $P_{\pm}$, and we say the $C_i$ form a standard triple of mutually tangent circles.  Otherwise, the circles $C_1, C_2, C_3$ are a \emph{non-standard triple} of mutually tangent circles (see Figure \ref{fig:NSMT} parts $(a)$ and $(b)$ for two views). 

Non-standard triples admit a convenient combinatorial description as well.  The three points of tangency in any set of three mutually tangent circles correspond to three edges $e_1,e_2,e_3$ in the crushtacean $\Gamma^{\ast}$ which form a 3-edge cut.  If the set is a non-standard triple, then the three edges are not adjacent, and neither component of $\Gamma^{\ast}\setminus \{e_1,e_2,e_3\}$ is a single vertex.  A \emph{non-trivial 3-edge cut} of $\Gamma^{\ast}$ will be one in which neither component of $\Gamma^{\ast}\setminus \{e_1,e_2,e_3\}$ is a single vertex (see Figure \ref{fig:NSMT}$(c)$).  Thus every non-standard triple in the circle packing corresponds to a non-trivial 3-edge cut of the crushtacean, and observe that the converse is also true.

Note also that non-standard three-edge cuts have an odd number of painted edges.  Trivalent graphs have an even number of vertices, since the sum of the degrees is both $3V$ and $2E$ (where $V, E$ represent the number of vertices and edges, respectively).  Given a three-edge cut, there are an odd number of vertices in each component of the cut graph since replacing a component with a single vertex results in another trivalent graph with an even number of vertices.

Two upcoming results are easily stated in terms of non-trivial 3-edge cuts. First, Corollary \ref{cor:DisCDCircPackCrush} shows that non-trivial 3-edge cuts in $\Gamma^{\ast}$ correspond to one class of thrice-punctured spheres in FAL complements.  Moreover, Theorem \ref{thm:CrushCharPrime} shows that an FAL complement is belted-sum prime if and only if every non-trivial 3-edge cut of the crushtacean is fully painted (i.e. all edges in the cut are painted).

\section{Thrice-punctured spheres}\label{sec:3ps}

This section classifies embedded, totally geodesic, thrice-punctured spheres in FAL complements. Adams proves that every essential, embedded thrice-punctured sphere in an orientable hyperbolic 3-manifold is totally geodesic (see \cite{ad1}), so assume all thrice-punctured spheres in this work are such. We begin by specifying two closely related classes of them, non-standard and non-reflection disks.  

Recall that, by convention, an FAL comes with a specific choice of reflection surface and a choice of crossing disk for each crossing circle relative to it.  Such choices lead to the standard cell decomposition $\mathcal{C}$ and standard polyhedra $P_{\pm}$ of Section \ref{sec:FAL}. Thrice-punctured spheres in $\mathcal{C}$ will be called \emph{standard disks}. The chosen crossing disks are standard, as are any thrice-punctured sphere components of the reflection surface.  These are special since they form faces on $P_{\pm}$. \emph{Non-standard disks} are thrice-punctured spheres not in $\mathcal{C}$.  The second class of disks, called \emph{Non-reflection disks}, are thrice-punctured spheres that are not part of the reflection surface of an FAL complement.  Thus non-reflection disks are the non-standard disks together with the chosen crossing disks.

The main result of this section is Theorem \ref{thm:NR3ps}, which combines results of this section into a classification of non-reflection disks in an FAL complement.  The classification is based on two properties: the intersection of non-reflection disks with the reflection surface and the slopes of their punctures along components of the FAL. We begin by describing possible intersections with the reflection surface.

\begin{center}
\begin{figure}[h]
\[
\begin{array}{ccc}
\includegraphics[width=1.75in]{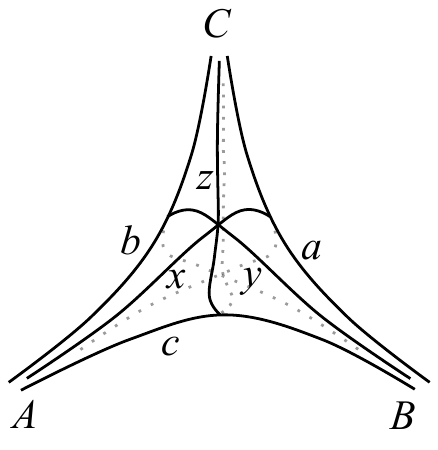}&\hspace{0.5in}&\includegraphics[width=1.75in]{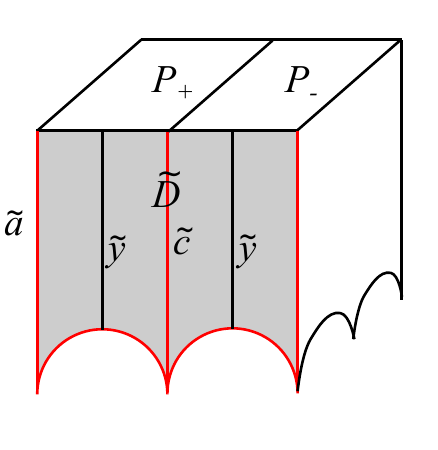}\\
(a)\textrm{ Geodesics on }S & &(b)\textrm{ Standard domain }\mathcal{F} 
\end{array}
\]
\caption{Neighborhood of $C$ in $P_+$}
\label{fig:3psGeod}
\end{figure}
\end{center}

Both the reflection surface and a non-reflection disk are embedded, totally geodesic surfaces.  If two such surfaces intersect non-trivially, they do so in a union of pairwise disjoint, simple geodesics on each. Thrice-punctured spheres have only the six simple geodesics pictured in Figure \ref{fig:3psGeod}$(a)$.  Those labeled $a, b, c$ are \emph{non-separating} geodesics since slicing along one does not separate the sphere. Slicing along $x$, $y$ or $z$ separates $S$, so they are \emph{separating} geodesics. We shall see that, if $\mathcal{A}_L$ is any FAL other than the fully-twisted Borromean rings, then a non-reflection disk $S$ is always orthogonal to $R$.  The geodesics in $S\cap R$ help distinguish non-reflection disks.   Propositions \ref{prop:2PiDisk}, \ref{prop:3psDisjointFromCD}, and \ref{prop:3psSeparatesD} will show that, if $\mathcal{A}_L$ has at least three crossing circles, then $S\cap R$ will either consist of all the non-separating geodesics on $S$ (termed $n$-disks), or exactly one separating geodesic of $S$ (called $s$-disks).

Each puncture of $S$ induces a slope on a torus neighborhood of one component of $\mathcal{A}_L$, which also helps classify non-reflection disks. Recall that a \emph{slope} is the isotopy class of an unoriented simple closed curve on a torus. Given a component $J$ of $\mathcal{A}_L$, let $T_J\subset \mathbb{S}^3$ denote the torus boundary of a tubular neighborhood $V(J)$ of $J$. Since $J$ is a knot in $\mathbb{S}^3$, there is a natural meridional slope $m$ and a longitudinal slope $\ell$ on $T_J$, determined by the following properties:  $m$ links $J$ once and intersects $\ell$ once, while $\ell$ is homologous to $J$ in $V(J)$ and null-homologous in $\mathbb{S}^3\setminus V(J)$ (see \cite[Theorem 3.1]{bz}).  We refer to any slope of the form $\ell \pm n\cdot m$ as a \emph{generalized longitude}.  Geometrically, $J$ corresponds to a cusp of $M$, and we think of $T_J$ as the boundary of a neighborhood of this cusp.  If an embedded totally geodesic surface $F$ is punctured by $J$, we frequently refer to the puncture by the slope $F\cap T_J$.  For example, a crossing disk of $M$ has two meridional knot circle punctures and one longitudinal crossing circle puncture.  Alternatively, consider a component $R_1$ of the reflection surface.   It has (generalized) longitudinal punctures along knot circles, and two types of punctures along a crossing circles, depending on whether the crossing circle is flat or twisted. If a flat crossing circle punctures $R_1$, it does so in meridian(s) of the crossing circle.  On the other hand, placing a half-twist in a crossing circle connects the two meridians, forming a single puncture of $R_1$ consisting of one longitude plus two meridians.  These observations will be useful in proving Lemma \ref{lem:Standard}. Types of punctures, then, also help classify non-reflection disks.

The main tool used for studying non-reflection disks will be their pre-images in a standard domain $\mathcal{F}$.  We begin by describing the preimages of chosen crossing disks, the simplest non-reflection disks.  Let $D$ be a crossing disk in an FAL complement, and consider how geodesics in a standard domain $\mathcal{F}$ project to the simple geodesics on $D$. Observe that the non-separating geodesics on $D$ are the one-cells of $\mathcal{C}$ since they join distinct punctures.  Thus the edges of $\mathcal{F}$ project to non-separating geodesics on crossing disks. Separating geodesics of $D$ intersect non-separating geodesics orthogonally, so altitudes of shaded triangles on $\partial\mathcal{F}$ project to form separating geodesics in $D$ (here an altitude of an ideal triangle is the geodesic ray orthogonal to one side and toward the opposite vertex). The altitudes labeled $\tilde{y}$ in Figure \ref{fig:3psGeod}$(b)$, for example, project to the separating geodesic $y$ of Figure \ref{fig:3psGeod}$(a)$.  The preimage of a crossing disk in $\mathcal{F}$, then, lies in the boundary of a standard domain and is well understood.

For other non-reflection disks we apply techniques of Knavel and Trapp in \cite{kt}, who studied embedded, totally geodesic surfaces in FAL complements by analyzing their pre-images in a standard domain $\mathcal{F}$. We review some relevant results from \cite{kt}, and begin with their definition of a geodesic disk. 

\begin{defn}\label{defn:GeoDisk}
Let $S\subset M$ be an embedded totally geodesic surface in the FAL complement $M$, and let $\widetilde{S}$ denote its preimage in the universal cover. A \emph{geodesic disk} is a connected component of $\widetilde{S}\cap P_+$, or of $\widetilde{S}\cap P_-$.  

Further, faces of the polyhedra $P_{\pm}$ are \emph{standard geodesic disks} since they project to the standard cell decomposition $\mathcal{C}$.  A \emph{non-standard geodesic disk} is one that is not a face of $P_{\pm}$. 
\end{defn}

\begin{figure}[h]
\[
\begin{array}{ccc}
\includegraphics[width=1.25in]{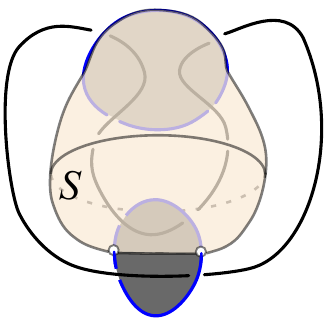}&\includegraphics[width=1.6in]{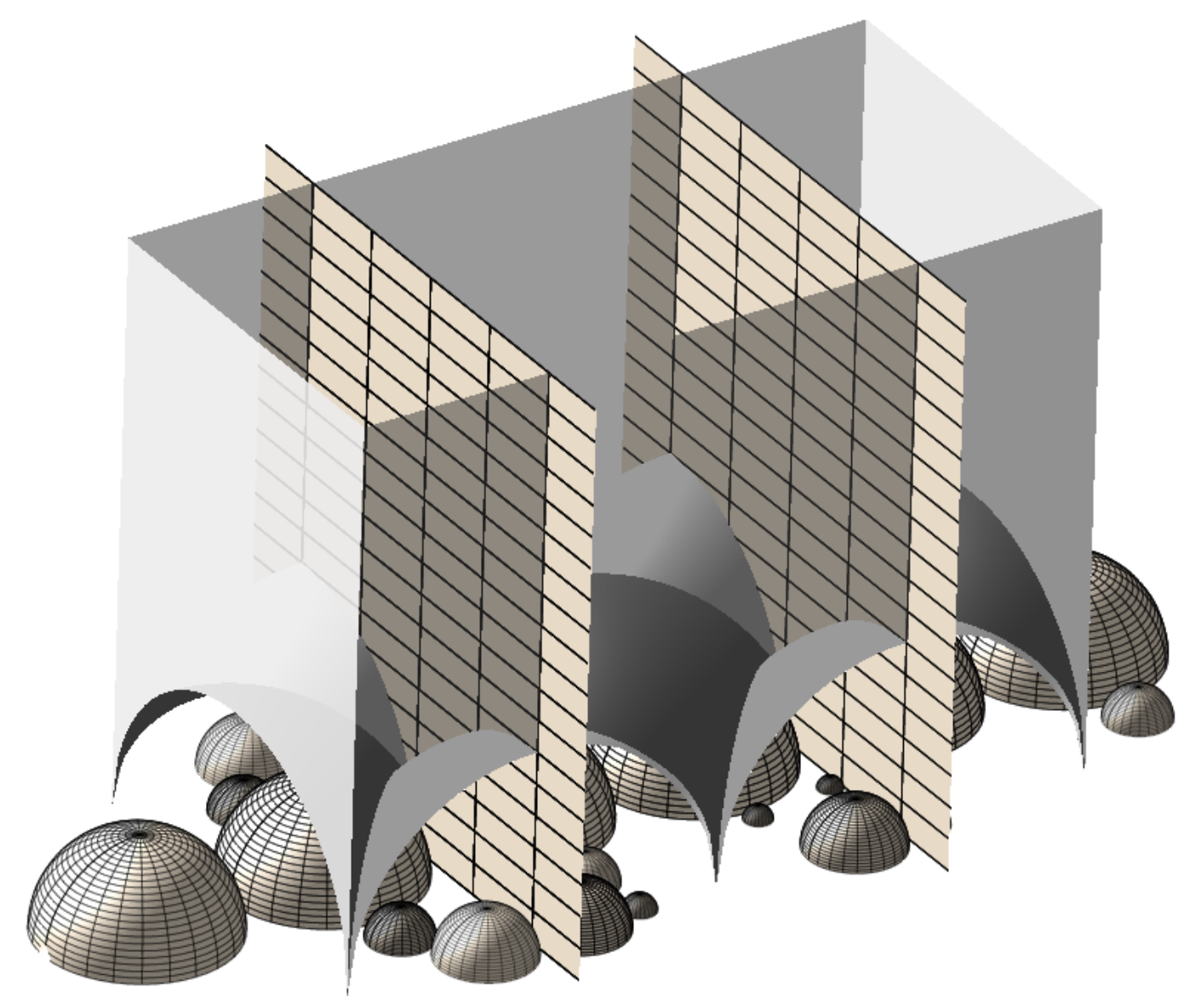}& \includegraphics[width=1.55in]{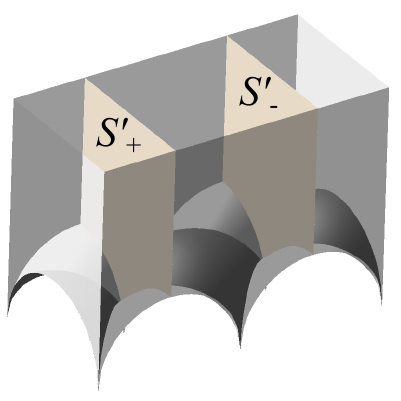}\\
(a)\textrm{ A non-reflection} & (b)\textrm{ Universal cover }& (c)\textrm{ Geodesic disks}\\
\phantom{.}\textrm{      thrice-punctured} &\textrm{ with preimage }\widetilde{S}\textrm{ and}&\widetilde{S}\cap\mathcal{F}\\
\textrm{ sphere} &\textrm{  standard domain }\mathcal{F} &
\end{array}
\]
\caption{Non-reflection sphere $S$, its $\widetilde{S}$, and geodesic disks}
\label{fig:GeoDisk}
\end{figure}

Since $S$ is an embedded totally geodesic surface, its preimage $\widetilde{S}$ in the universal cover is a union of disjoint planes in $\mathbb{H}^3$. Thus a geodesic disk is the intersection of a plane $\mathcal{P}$ in $\mathbb{H}^3$ with either $P_+$ or $P_-$. The boundary $\gamma$ of a geodesic disk $S'_+$ in $P_+$, say, is a planar polygon given by $\gamma = \mathcal{P}\cap\partial P_+$.   For example, Figure \ref{fig:GeoDisk}$(a)$ depicts a non-reflection disk in $M$, with some planes of its preimage $\widetilde{S}$ in the universal cover found in Figure \ref{fig:GeoDisk}$(b)$.  The geodesic disks $S'_{\pm}$ are highlighted in Figure \ref{fig:GeoDisk}$(c)$, both quadrilaterals with two ideal and two material vertices.

The boundary of a non-standard geodesic disk can intersect a shaded triangle $T$ of $\partial P_{\pm}$ only along an edge or altitude.  To see this, let $D$ be a crossing disk and recall that its preimage in $\partial P_{\pm}$ consists of two pairs of shaded triangles (see Figure \ref{fig:3psGeod}$(b)$).  Edges of $\partial P_{\pm}$ project to the non-separating, and altitudes of shaded triangles project to the separating, geodesics of $D$.  Now a non-reflection sphere $S$ is an embedded, totally geodesic surface, as are components of $\mathcal{C}$. Thus $S$ must intersect each component of $\mathcal{C}$ in a collection of disjoint, simple geodesics on each surface.  In particular, if $S\cap D$ is non-empty, then it is a collection of disjoint separating and non-separating geodesics on each surface.  But then $\mathcal{P}\cap T$ must be an edge or altitude of $T$.  

With these initial observations in hand, we classify the types of geodesic disks that result from non-reflection disks. Area considerations limit the types of geodesic disks that can arise from a non-standard thrice-punctured sphere.  Recall that the area of a hyperbolic polygon with external angles $\epsilon_i$ is $ \sum_{i=1}^n \epsilon_i - 2\pi$, with the convention that $\epsilon = \pi$ for ideal vertices. We begin by demonstrating the general observation that the area of a geodesic disk in an FAL fundamental domain is always an integer multiple of $\pi$.

\begin{lemma}\label{lem:AreaNPi}
Let $S'\subset P_{\pm}$ be a non-standard geodesic disk for the embedded totally geodesic surface $S$ in an FAL complement.  Then

\begin{enumerate}[label=\roman*.]
\item Material vertices of $\gamma=\partial S'$ are right-angled,
\item Material vertices can be paired so that their adjacent altitudes share an ideal vertex, and
\item The area of $S'$ is an integer multiple of $\pi$.
\end{enumerate} 
\end{lemma}

\begin{proof}
Let $S'\subset P_{\pm}$ be a non-standard geodesic disk for the embedded totally geodesic surface $S\subset M$, with boundary polygon $\gamma = S'\cap \partial P_{\pm}$. 

Every edge of $P_{\pm}$ bounds one unshaded face and one shaded triangle. If $q$ is a material vertex of $\gamma$, then one edge $\gamma_0\subset \gamma$ intersects the interior of a shaded triangular face $T_0\subset\partial P_{\pm}$ (see Figure \ref{fig:AreaNPi}). As above, the geodesic ray $\gamma_0$ is an altitude of $T_0$ and projects to (half of) a separating geodesic on some crossing disk $D\subset M$.  This is precisely the situation of \cite[Lemma 4.2.i.]{kt}, which shows that $S'$ is orthogonal to both faces of $\partial P_{\pm}$ incident with $q$ (labeled $T_0$ and $R$ in Figure \ref{fig:AreaNPi}). Then $q$ is the intersection of three pairwise orthogonal planes in $\mathbb{H}^3$, one of which contains $S'$, so the angle of $S'$ at $q$ is a right angle.

\begin{figure}[h]
\begin{center}
\includegraphics{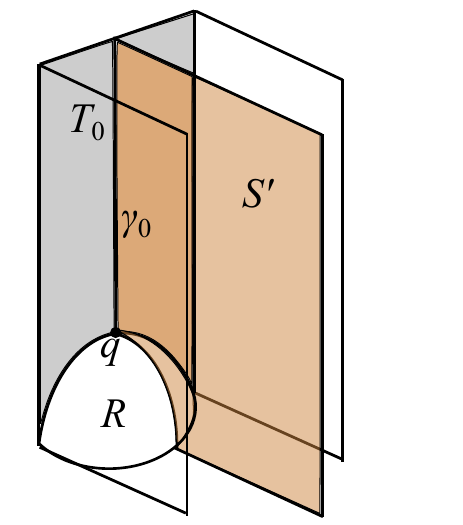}
\end{center}
\caption{Geodesic disk $S'$ with material vertex $q$}
\label{fig:AreaNPi}
\end{figure}

To see that material vertices come in pairs, let $p$ be the ideal ``endpoint" of $\gamma_0$, and let $T_1$ be the shaded triangle in $\partial P_{\pm}$ that shares $p$ with $T_0$. Now apply a M\'obius transformation to $\mathcal{F}$ so that $p$ is at infinity.  The geodesic disk $S'$ is orthogonal to $T_0$ and $P_{\pm}$ is right-angled, so $S'$ is orthogonal to the face $T_1$ opposite to $T_0$ at $p$.  Thus $S'\cap T_1$ is an altitude of $T_1$ with one material endpoint $q'$, and there are an even number of material endpoints on $\gamma$.   Further, the altitudes with base material vertices $q,q'$ share the ideal vertex $p$ on $\gamma$.

Since each ideal vertex of $\gamma$ contributes $\pi$ to the area of $S'$ and there are an even number of material vertices, each contributing $\pi/2$, the area of $S'$ is an integer multiple of $\pi$. 
\end{proof}

Lemma \ref{lem:AreaNPi} has the following immediate corollary.

\begin{cor}\label{cor:PossGeodDisk}
If $S$ is an embedded, totally geodesic thrice-punctured sphere in an FAL complement, then $S$ has either two geodesic disks of area $\pi$ or one of area $2\pi$.
\end{cor}

\begin{proof}
Let $S$ is an embedded, totally geodesic thrice-punctured sphere in an FAL complement.  Then $S$ has area $2\pi$, so Lemma \ref{lem:AreaNPi} $(iii)$ implies it either has two geodesic disks of area $\pi$ or one of area $2\pi$.
\end{proof}

The first proposition shows that there is a unique FAL complement that contains a non-reflection disk with a single area $2\pi$ geodesic disk.  

\begin{prop}\label{prop:2PiDisk}
If $S$ is a non-standard disk in the FAL complement $M$ with a single geodesic disk of area $2\pi$, then $M$ is the fully twisted Borromean rings.
\end{prop}

\begin{proof}
Suppose that $S'$ is a non-standard geodesic disk with area $2\pi$ that projects to a non-reflection disk $S$.  Up to a reflection of $M$ we can suppose $S'\subset P_+$.  If $q$ were a material vertex of $\partial S'$, the unshaded face adjacent to $q$ contains an edge of $\partial S'$ which glues $S'$ to another geodesic disk.  Since $S'$ is the only geodesic disk for $S$, this implies all of its boundary vertices must be ideal. Thus $\gamma = \partial S'$ is an ideal square (see Figure \ref{fig:2Pi3ps}$(a)$).

\begin{center}
\begin{figure}[h]
\[
\begin{array}{ccc}
 \includegraphics[width=1.3in]{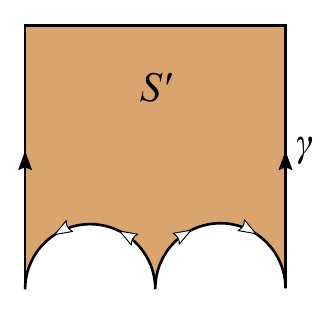}&\includegraphics[width=1.2in]{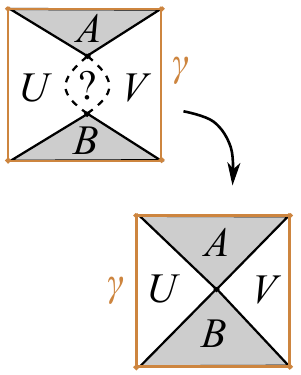}&\includegraphics[width=1.2in]{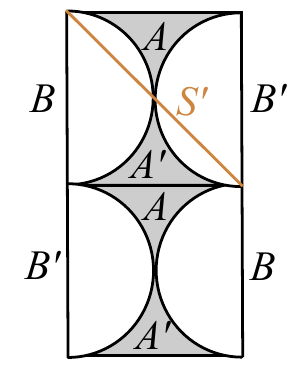} \\
(a)\textrm{ Gluing on }S'&(b) \textrm{ One side of }\gamma & (c)\ S'\textrm{ in }\mathcal{F}
\end{array}
\]
\caption{Area $2\pi$ geodesic disk}
\label{fig:2Pi3ps}
\end{figure}
\end{center}

The edges of $\gamma$ lie on distinct faces since the plane containing $S'$ intersects a face of $\partial P_+$ in at most one geodesic arc.  In addition, vertices of $P_+$ are 4-valent and checkerboard colored, so the faces of $\partial P_+$ on one side of $\gamma$ alternate between shaded triangles and unshaded faces as $\gamma$ is traversed. Let $A$ and $B$ be the shaded triangles of $P_+$ adjacent to one side of $\gamma$, and let $U, V$ denote the unshaded faces adjacent to the same side (see Figure \ref{fig:2Pi3ps}$(b)$).  Both $U$ and $V$ share an edge with the shaded triangle $A$, so they share the vertex of $A$ not on $\gamma$.  Similarly, they share the third vertex of $B$. Now $U$ and $V$ share at most one ideal vertex (since their corresponding circles in the circle packing are tangent at most once), so the third vertices of shaded triangles $A$ and $B$ coincide.  Thus one side of $\gamma$ consists solely of the four ideal triangles $A, U,B,V$.

The same argument applies to the other side of $\gamma$, making $P_{\pm}$ regular ideal octahedra. Then $\mathcal{F}$ is the standard domain for the Borromean rings or its half-twist partners. Moreover, $\mathcal{F}$ can be chosen so that $S'$ is on the vertical half plane that intersects $P_+$ along a diagonal as in Figure \ref{fig:2Pi3ps}$(c)$. 

To see that the gluing pattern on $\mathcal{F}$ is that of the fully twisted Borromean rings (see Figure \ref{fig:F}$(c)$), apply a M\"obius transformation to $\mathcal{F}$ so that the two vertical edges of $S'\subset P_+$ are identified as in Figure \ref{fig:2Pi3ps}$(a)$.  The vertical shaded triangles on the left of $\mathcal{F}$ are glued to ones on the right by parabolic isometries that preserve infinity.  Since the gluing maps also identify the vertical edges of $\gamma$, the faces labeled $B$ in Figure \ref{fig:2Pi3ps}$(c)$ must be identified, resulting in a twisted crossing disk.  A similar argument shows the gluing pattern on $A$ faces results in a twisted crossing disk as well, finishing the proof. 
\end{proof}

Ideal square geodesic disks occur for more general embedded, totally geodesic surfaces in FAL complements.  The point of Proposition \ref{prop:2PiDisk}, however, is that in most cases they are glued to other geodesic disks, giving the resulting surfaces too much area to be a thrice-punctured sphere.  See \cite[Example 5.2]{kt} for examples of surfaces with ideal square geodesic disks in flat, 2-bridge FAL complements.  A more thorough discussion of thrice-punctured spheres in the fully-twisted Borromean rings complement precedes Lemma \ref{lem:Not2Pi3ps}.

Now consider area $\pi$ geodesic disks.  An area $\pi$ geodesic disk must have exterior angle sum $3\pi$, and it either has material vertices or not.  Lemma \ref{lem:AreaNPi} showed that material vertices are right-angled and come in pairs, with an ideal vertex between them, thus contributing $2\pi$ to the angle sum.  The area of a geodesic disk with more than two material vertices, then, is at least $3\pi$ and cannot project to a thrice-punctured sphere.  Thus an area $\pi$ geodesic disk that projects to a non-standard disk will either be an ideal triangle (with no material vertices), or a quadrilateral whose boundary contains two material vertices separated by two ideal vertices.  The term \emph{material-vertex disk} refers to an area $\pi$ geodesic disk with two material and two ideal vertices (see Figure \ref{fig:GeoDisk}$(c)$ for an example of material-vertex disks), while a \emph{non-standard ideal triangle} is an ideal triangle geodesic disk that is not a face of $P_{\pm}$

\begin{figure}[h]
\[
\begin{array}{ccc}
\includegraphics[width=1.5in]{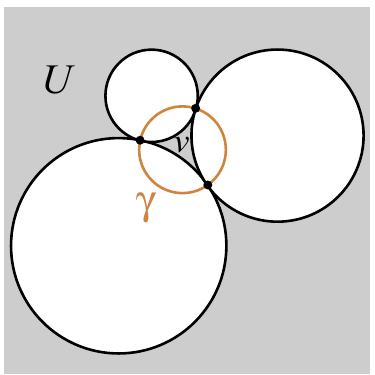}&\hspace{0.5in}& \includegraphics[width=1.5in]{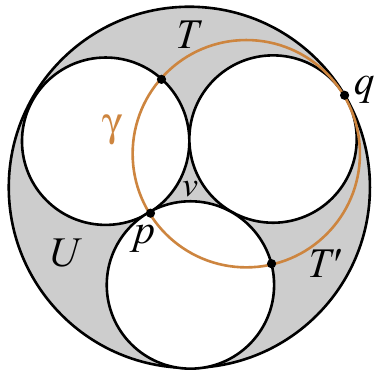}\\
(a) \textrm{ An ideal triangle} && (b) \textrm{ A material-vertex disk}
\end{array}
\]
\caption{Area $\pi$ geodesic disks in the circle packing}
\label{fig:AreaPiDisks}
\end{figure}

A non-standard geodesic disk results from intersecting $P_{\pm}$ with a plane $\mathcal{P}$ in $\mathbb{H}^3$, where $\mathcal{P}$ projects to a non-standard disk in an FAL complement.  We describe area $\pi$ disks in terms of the circle packing for $P_{\pm}$ and the boundary circle for $\mathcal{P}$ (see Figure \ref{fig:AreaPiDisks}).  The circle packing for $P_{\pm}$ must have three mutually tangent circles for a non-standard ideal triangle to exist. Moreover, the interstices (regions $U$ and $v$ in Figure \ref{fig:AreaPiDisks}$(a)$) between the circles must contain non-trivial portions of the circle packing, otherwise the ideal triangle $\mathcal{P}\cap P_{\pm}$ is a face of $P_{\pm}$. Thus non-standard ideal triangles correspond to the non-standard triples of Section \ref{sec:FAL}. The four circles involved in a material-vertex disk (two at each ideal point of its boundary) are mutually tangent in the circle packing, forming a Descartes quadruple.  The interstices $U$ and $v$ of Figure \ref{fig:AreaPiDisks}$(b)$ may contain non-trivial portions of the circle packing, while $T$ and $T'$ are shaded triangular faces of $P_{\pm}$ containing the altitudes in the boundary of the disk.

The above discussion is summarized in the following characterization of area $\pi$ non-standard geodesic disks.

\begin{lemma}\label{lem:AreaPiDisks}
Let $M$ be an FAL complement with standard polyhedra $P_{\pm}$, and let $S'$ be an area $\pi$ geodesic disk for a non-standard disk $S\subset M$.  

Then $S'$ is orthogonal to $\partial P_{\pm}$ and is either a non-standard ideal triangle or a material-vertex disk.
\end{lemma}

\begin{proof}
The preceding discussion restricts area $\pi$ disks to non-standard ideal triangles or material-vertex disks, so it must be shown that
both types are orthogonal to the standard cell decomposition.  The discussion leading up to the definition of non-standard triples shows that a non-standard ideal triangle meets $\mathcal{C}$ orthogonally.  Further, the boundary of a material-vertex disk contains two altitudes of shaded triangles, and \cite[Lemma 4.2.i.]{kt} states that any geodesic disk whose boundary contains an altitude meets $\mathcal{C}$ orthogonally. 
\end{proof}

\begin{figure}[h]
\[
\begin{array}{ccc}
\includegraphics[width=2in]{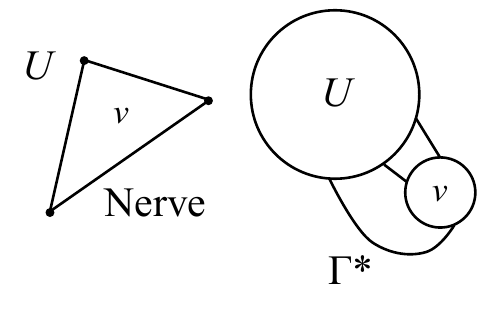}&\hspace{0.25in}& \includegraphics[width=2in]{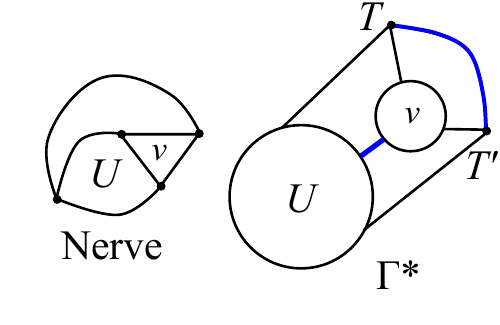}\\
(a) \textrm{ An ideal triangle} && (b) \textrm{ Material-vertex disk}
\end{array}
\]
\caption{Recognizing disks from graphs}
\label{fig:DisksFromGraphs}
\end{figure}

Area $\pi$ disks can also be recognized in the nerve and crushtacean of an FAL.  Indeed, we have already seen that a non-standard triple corresponds to a non-trivial 3-edge cut in the crushtacean.  A material-vertex disk corresponds to a $K_4$ subgraph in the nerve, two of whose regions are faces of the triangulation while the other two regions may contain additional faces in the triangulation.  This corresponds to two 3-edge cuts in the crushtacean that share an edge (the edge between $U$ and $v$ in Figure \ref{fig:DisksFromGraphs}$(b)$).  

The following proposition characterizes non-standard disks whose geodesic disks are ideal triangles.  It is interesting to note that every non-standard ideal triangle and its reflection glue to form an non-standard disk.  In contrast, whether or not material-vertex disks project to non-standard disks depends on how faces of standard polyhedra are identified. 

\begin{prop}\label{prop:3psDisjointFromCD}
Let $M=\mathbb{S}^3\setminus\mathcal{A}$ be an FAL complement with reflection surface $R$ and standard domain $\mathcal{F}=P_+\cup P_-$.  

The surface $S\subset M$ is a non-standard disk disjoint from chosen crossing disks if and only if it is the projection of a non-standard ideal triangle and its reflection in $P_{\pm}$.  

Moreover, in this case $S$ satisfies:
\begin{enumerate}[label=\roman*.]
\item $S$ intersects $R$ orthogonally in its three non-separating geodesics.
\item The punctures of $S$ are either one crossing-circle longitude and two knot-circle meridians, or three longitudes on distinct crossing circles.
\end{enumerate}
\end{prop}

\begin{proof}
Let $M$ be an FAL complement and $S$ a non-standard disk in $M$ that is disjoint from the crossing disks. Non-standard ideal triangles are disjoint from shaded triangular faces of $P_{\pm}$, so their projections are disjoint from chosen crossing disks in $M$. Since area $2\pi$ geodesic disks and material-vertex disks always intersect shaded faces of $P_{\pm}$, their projections intersect chosen crossing disks (by Proposition \ref{prop:2PiDisk} and the definition of material-vertex disks). Thus the geodesic disks for $S$ must be non-standard ideal triangles.   

Conversely we must show that each non-standard ideal triangle and its reflection, denoted $S_{\pm}'$, project to a non-standard disk $S\subset M$. By Lemma \ref{lem:AreaPiDisks} the disks $S_{\pm}'$ are orthogonal to $\partial P_{\pm}$, so they project to an embedded, totally geodesic surface $S\subset M$ which is disjoint from chosen crossing disks.  Moreover, the surface $S$ is two ideal triangles with edges glued to their reflections, resulting in a thrice-punctured sphere and proving the first claim of the proposition.

We now verify the stated properties of $S$.  Note that slicing the projection of $S'_+\cup S'_-$ along one geodesic of $\partial S'_{\pm}$ does not separate it; therefore $S$ intersects the reflection surface along its non-separating geodesics.  Moreover, we've already noted that $S'_{\pm}$ are orthogonal to $\partial P_{\pm}$, so $S$ is orthogonal to the reflection surface.

We now address the punctures of $S$.  Since $S$ is disjoint from all crossing disks, it will have the same meridional and longitudinal punctures as a crossing disk.  Hence knot circle punctures of $S$ will be meridians while crossing circle punctures are longitudinal.  

To see that crossing circle punctures of $S$ must be on distinct cusps note that exactly one vertex of $P_+$ projects to each crossing circle cusp of $M$ (and similarly for $P_-$).  Thus vertices of $S_{\pm}'$ that project to crossing circle cusps correspond to distinct crossing circles as desired.

It remains to show that $S$ has an even number of knot circle punctures. If all punctures of $S$ were meridians of knot circles, then $S$ would be a sphere in $S^3$ punctured three times by a link, which is impossible.  Similarly, if $S$ had exactly one knot circle puncture, then the other two punctures would be longitudes of crossing circles.  Gluing the two crossing disks to $S$ would result in a sphere in $S^3$ punctured five times by a link, arriving at the same contradiction.  Thus $S$ has an even number of knot circle punctures, completing the proof.
\end{proof}

Proposition \ref{prop:3psDisjointFromCD} gives a very precise description of non-standard disks which are disjoint from crossing disks.  Notice that, like crossing disks, they intersect the reflection surface orthogonally in their non-separating geodesics.  For this reason we refer to them as $n$-disks, the ``$n$" emphasizing this observation. There are two kinds of $n$-disks though, depending on their punctures, and we make the following formal definition.  

\begin{center}
\begin{figure}[h]
\[
\begin{array}{ccc}
\includegraphics{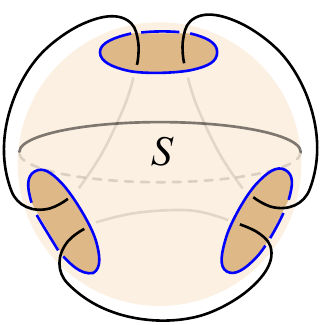}&&\includegraphics{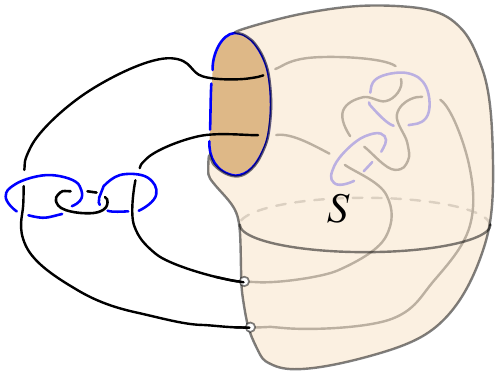}\\
(a)\textrm{ Longitudinal disk}&&(b)\textrm{ (Generalized) crossing disk} \\
\end{array}
\]
\caption{Two types of $n$-disks}
\label{fig:NDisks}
\end{figure}
\end{center}

\begin{defn}\label{defn:NDisks}
Let $M$ be an FAL complement with chosen reflection surface $R$. An $n$-disk is a non-reflection disk $S\subset M$ is one that intersects $R$ orthogonally along its three non-separating geodesics.  There are two types of $n$-disks:
\begin{enumerate}[label=\roman*.]
\item A \emph{longitudinal disk} is an $n$-disk that has three longitudinal punctures on distinct crossing circles (see Figure \ref{fig:NDisks}$(a)$).
\item A \emph{crossing disk} is an $n$-disk that has one crossing-circle longitudinal and two knot-circle meridional punctures (see Figure \ref{fig:NDisks}$(b)$).
\end{enumerate} 
\end{defn}

Note that Definition \ref{defn:NDisks}$(ii)$ generalizes the term ``crossing disk".  In fact, generalized crossing disks are potential crossing disks for the crossing circle.  Generalized crossing disks will play a crucial role in classifying belted sum decompositions of FALs, as well as in defining a canonical one.

The fact that non-standard ideal triangles can be recognized in circle packings and crushtaceans leads immediately to the following corollary (see Figures \ref{fig:AreaPiDisks}$(a)$ and \ref{fig:DisksFromGraphs}$(a)$).  

\begin{cor}\label{cor:DisCDCircPackCrush}
Let $M=\mathbb{S}^3\setminus\mathcal{A}$ be an FAL complement.  Non-standard disks in $M$ which are disjoint from crossing disks are in one-to-one correspondence both with non-standard triples in the circle packing and with non-trivial 3-edge cuts of the crushtacean.

Furthermore, (non-standard) crossing disks correspond to non-trivial 3-edge cuts with one painted edge while longitudinal disks correspond to those with all edges painted.
\end{cor}

\begin{proof}
The follows from simply interpreting Proposition \ref{prop:3psDisjointFromCD} in terms of the circle packing and crushtacean.  Do recall that painted edges in the crushtacean represent vertices of $P_{\pm}$ that correspond to crossing circles, so the number of painted edges in the 3-edge cut will correspond to the number of crossing circle punctures on the disk.
\end{proof}

\begin{figure}[h]
\begin{center}
\includegraphics{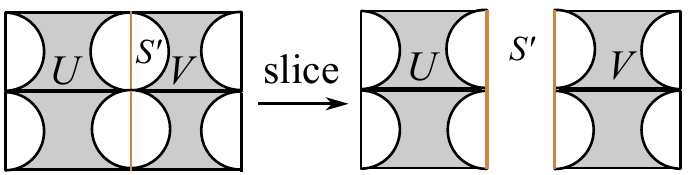}
\end{center}
\caption{Slicing a standard domain along an $n$-disk}
\label{fig:IdealTriFD}
\end{figure}

One slices along thrice-punctured spheres when performing belted-sum decompositions in FAL complements.  In a standard domain this corresponds to slicing along geodesic disks.  Therefore, it will be helpful to analyze how non-standard ideal triangles and their reflections separate standard domains.  Figure \ref{fig:IdealTriFD} illustrates a standard domain $\mathcal{F}=P_+\cup P_-$ containing the union $S'$ of a non-standard ideal triangle and its reflection.  The shaded regions labeled $U$ and $V$ contain additional circles in the circle packing.  We can assume that the vertex at infinity corresponds to a crossing circle cusp, since either one or all vertices of $S'$ correspond to crossing circle cusps by Proposition \ref{prop:3psDisjointFromCD}.

We now analyze how the gluing instructions on $\mathcal{F}$ identify faces of the pieces after slicing along $S'$.  Unshaded faces of $P_{\pm}$ glue to their reflections, so unshaded faces glue within pieces.  Shaded triangles are paired with adjacent ones to determine gluings.  Shaded triangles which share a vertex on $\partial S'$, then are the only ones which can identify faces between the pieces.  In the case of a longitudinal disk, then, all shaded triangles sharing a vertex of $S'$ are paired, providing many identifications between the pieces of $\mathcal{F}\setminus S'$.  On the other hand, if $S'$ projects to a generalized crossing disk only infinity corresponds to a crossing circle and only the vertical shaded faces of $\mathcal{F}$ glue the two pieces together.

Having characterized non-reflection thrice-punctured spheres with ideal triangle geodesic disks, we turn our attention to material-vertex disks.  Material-vertex disks have two altitudes in their boundaries which share an ideal vertex.  The next proposition shows that a material-vertex disk (together with its reflection) projects to a non-standard disk if and only if the ideal vertex shared by these altitudes corresponds to a flat crossing circle. In this case, the remaining ideal vertex corresponds to a crossing circle as well.

\begin{prop}\label{prop:3psSeparatesD}
Let $M=\mathbb{S}^3\setminus\mathcal{A}$ be an FAL complement with standard domain $\mathcal{F}=P_+\cup P_-$.  

The surface $S\subset M$ is a non-reflection thrice-punctured sphere that separates a crossing disk $D$ if and only if $S$ is the projection of a material-vertex disk and its reflection in $P_{\pm}$ whose altitude edges project to a separating geodesic on the flat crossing disk $D$.

In this case, $S$ intersects the reflection surface orthogonally along one of its separating geodesics.  Moreover, the punctures of $S$ are a longitude of one crossing-circle and two meridians of a different crossing-circle.
\end{prop}

\begin{proof}
Let  $S\subset M$ be a non-reflection thrice-punctured sphere that intersects the crossing disk $D$ along one of its separating geodesics $\gamma_D$.  Since the preimage in $P_{\pm}$ of $\gamma_D$ consists of altitudes of shaded triangles, at least one geodesic disk has material vertices, call it $S'_+$. Then $S'_+$ is a material-vertex disk, and must glue to its reflection $S'_-$ along the rays interior to unshaded faces in $P_{\pm}$.  Area considerations imply these are the only two geodesic disks for $S$, so the altitudes must identify as well.  If $D$ is a twisted crossing disk, the identification of $S'_{\pm}$ results in a non-orientable surface.  Therefore $D$ is a flat crossing disk.

Conversely, suppose the material vertex disks $S'_{\pm}$ are reflections of each other, and their altitudinal boundary rays project to a separating geodesic on a flat crossing disk $D$.  This implies the altitudes of $\partial S'_+$ are identified, and the ideal vertex they share corresponds to the crossing circle $C$ which bounds $D$ (and similarly for the altitudes of $\partial S'_-$). The punctures of $S'_{\pm}$ shared by altitudes, then, project to two distinct meridional punctures of the crossing circle $C$. As in the proof of Lemma \ref{lem:AreaNPi}, applying \cite[Lemma 4.2.i.]{kt} we see that $S'_{\pm}$ is orthogonal to all faces of $\partial P_{\pm}$ that they meet.  Identifying edges, then, results in an embedded, totally geodesic surface $S$ in $M$. Moreover, the gluing pattern (altitudes glue within and other rays glue between geodesic disks) is that of a thrice-punctured sphere.  Note that first gluing altitudes within each $S'_{\pm}$ does not connect the disks.  Thus the unshaded boundary rays project to a separating geodesic of $S$, which is its intersection with the reflection surface of $M$.

It remains to show that the final puncture is a longitude of a crossing circle.  The remaining ideal vertex of $S'_+$ joins two geodesic rays interior to unshaded faces.  Therefore it and it's reflection identify to either a meridian of a knot circle or a longitude of a crossing circle.  If it were the meridian of a knot circle, then $S$ would be a sphere in $\mathbb{S}^3$ punctured three times by the link $\mathcal{A}$ (twice by $C$ and once by the knot circle), which is impossible.  Therefore, the remaining puncture is the longitude of a crossing circle.
\end{proof}

\begin{center}
\begin{figure}[h]
\[
\begin{array}{ccc}
\includegraphics{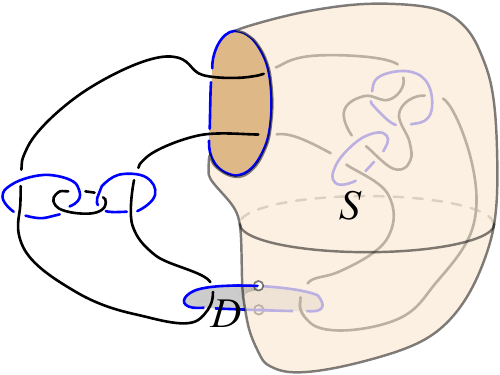}&\hspace{0.25in}&\includegraphics{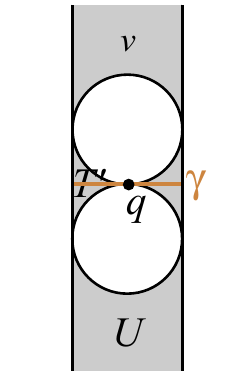}\\
(a)  \textrm{ Singly-separated disk}&&(b) \textrm{ $p$ at infinity}
\end{array}
\]
\caption{A singly-separated disk with circle-packing view}
\label{fig:SingSepDisk}
\end{figure}
\end{center}

The situation of Proposition \ref{prop:3psSeparatesD} is illustrated both in Figures \ref{fig:GeoDisk}$(a)$ and \ref{fig:SingSepDisk}$(a)$ and, as they appear frequently in the remainder of the paper, we make the following definition.

\begin{defn}\label{defn:SDisk}
Let $M$ be an FAL complement with chosen reflection surface $R$. A \emph{singly-separated disk}, or $s$-disk, is a non-reflection disk $S\subset M$ satisfying:
\begin{enumerate}[label=\roman*.]
\item $S$ intersects $R$ orthogonally,  
\item $S\cap R$ is one separating geodesic on $S$, and
\item $S$ has two meridional punctures along one crossing circle and one longitudinal puncture along a different crossing circle.
\end{enumerate} 
\end{defn}

Figures \ref{fig:AreaPiDisks}$(b)$ and \ref{fig:DisksFromGraphs}$(b)$ illustrated recognizing material-vertex disks from the circle packing and crushtacean.  The proof of Proposition \ref{prop:3psSeparatesD} justifies the edge-paintings of Figure \ref{fig:DisksFromGraphs}$(b)$, showing that the vertices corresponding to painted edges do correspond to crossing circles.  Transferring vertex $p$ of Figure \ref{fig:AreaPiDisks}$(b)$ to infinity via a M\"obius transformation yields Figure \ref{fig:SingSepDisk}$(b)$, from which it is easier to see that slicing along an $s$-disk separates $P_{\pm}$ into two pieces.  Moreover, the only gluing of faces between the two pieces of $P_+\setminus S_+'$ is along the vertical shaded triangular faces.  This is because the reflection between $P_{\pm}$ preserves the pieces of $P_{\pm}\setminus S_{\pm}'$ and unshaded faces glue to their reflections.  Further, shaded triangular faces are identified in adjacent pairs and, since the shaded faces intersected by $S_{\pm}'$ project to a flat crossing disk, the only pair separated by slicing along $S_{\pm}'$ is the vertical triangular faces. 

Combining Propositions \ref{prop:2PiDisk}, \ref{prop:3psDisjointFromCD}, and \ref{prop:3psSeparatesD} classifies non-reflection thrice-punctured spheres in FAL complements, proving the following theorem.

\begin{thm}\label{thm:NR3ps}
Let $S$ be a non-reflection thrice-punctured sphere in an FAL complement $M = \mathbb{S}^3\setminus \mathcal{A}$.  If $\mathcal{A}$ is not the fully twisted Borromean rings, then $S$ is orthogonal to the reflection surface in $M$ and is either a crossing, longitudinal or singly-separated disk.  
\end{thm}

\section{Separating Pairs}\label{sec:SP}

Now that we can recognize thrice-punctured spheres in FAL complements, we consider which spheres can be used to decompose FALs.  
Adams introduced a cut-and-paste method for constructing a hyperbolic manifold from two others containing incompressible, embedded thrice-punctured spheres (see \cite[Theorem 4.5]{ad1}). We review his construction here. Let $S_1\subset M_1$ and $S_2\subset M_2$ be incompressible, embedded thrice-punctured spheres in orientable hyperbolic three manifolds $M_1,M_2$.  Let $M_i' = M_i\setminus N(S_i)$ be the manifolds $M_i$ with a neighborhood of the spheres $S_i$ removed.  Thus the boundary of $M_i'$ consists of two thrice-punctured spheres which we label $S_i^0$ and $S_i^1$.  For $j \in \{0,1\}$, choose homeomorphisms $\lambda_j: S_1^j\to S_2^j$ which either both preserve or both reverse orientation.  Let $M$ be the manifold obtained by identifying the boundaries of $M_i'$ using the homeomorphisms $\lambda_j$. Theorem 4.5 of \cite{ad1} proves that the manifold $M$ constructed in this way is hyperbolic, and the volume of $M$ is the sum of the volumes of $M_1$ and $M_2$.  

As an application of this ``sum" of two manifolds, Adams introduces the notion of a belted sum of two links in $S^3$.  His definition of a belted sum is diagrammatic, saying that the link $L$ is the belted sum of links $L_1,L_2$ if they admit diagrams as in Figure \ref{fig:BS}.  We will use the notation $M=M_1\#_b M_2$ to denote a belted sum decomposition, where $M_i = \mathbb{S}^3\setminus L_i$ represent the link complements.

Clearly belted sums are special cases of punctured sphere sums in which the $S_i$ look like crossing disks in FALs in that they contain one longitudinal and two meridional punctures.  The homeomorphisms that yield the link $L$ identify punctures of the same type (meridional or longitudinal) without introducing additional crossings.  

We are concerned with this process in reverse and, in particular, applying it to fully augmented link complements.  It turns out that every punctured sphere decomposition of an FAL complement can be realized as a belted sum decomposition since the thrice-punctured spheres involved will be shown to always have one longitudinal and two meridional punctures. For this reason we ignore the more general setting and focus on belted-sum decompositions. By definition, if $M$ is a belted-sum of two manifolds then it contains a pair of disjoint thrice-punctured spheres that separate it.  To search for such decompositions, then, we introduce the notion of a separating pair of thrice-punctured spheres.

\begin{defn}\label{defn:SepPair}
Let $M$ be a hyperbolic FAL complement.  The pair $\{S_1,S_2\}$ of disjoint, essential, thrice-punctured spheres in $M$ is a \emph{separating pair} if $M\setminus\left(S_1\cup S_2\right)$ is disconnected.
\end{defn} 

We also formally define primality in the context of belted sums, and introduce notation.

\begin{defn}\label{defn:BSPrime}
A FAL is \emph{$b$-prime} if its complement is not the belted sum of two manifolds.
\end{defn}

Section \ref{sec:FALPrime} will use the absence of separating pairs to provide several characterizations of $b$-prime FALs, and Section \ref{sec:CanonDec} introduces a canonical belted sum decompositions of FALs. First, a clear understanding of separating pairs in FALs is required.

In the last section we showed that there are several types of thrice-punctured spheres in FAL complements: standard thrice-punctured spheres, the non-standard ones of Proposition \ref{prop:2PiDisk} that exist only in the fully twisted Borromean rings, and non-standard crossing disks, longitudinal disks, singly-separated disks (see Theorem \ref{thm:NR3ps}). Our goal is to determine when two thrice-punctured spheres form a separating pair.  The following lemma will prove a useful observation.

\begin{lemma}\label{lem:NotSP}
Let $\mathcal{A}$ be an FAL and let $\{S_1,S_2\}$ be a separating pair in $M = S^3\setminus \mathcal{A}$.  If the component $J\in\mathcal{A}$ punctures the pair $\{S_1,S_2\}$, then $J$ punctures the pair more than once.  
\end{lemma}
\begin{proof} 
Suppose that $\{S_1,S_2\}$ is a separating pair and that the component $J\in\mathcal{A}$ punctures $\{S_1,S_2\}$ exactly once, say on $S_1$.  Then the cusp of $M$ corresponding to $J$ has a neighborhood which is disjoint from $S_2$, and we let $T_J$ be the torus boundary of this cusp neighborhood. The curve  $S_1\cap T_J$ is a torus knot on $T_J$ whose complement on $T_J$ is an annulus.  This annulus is disjoint from $S_2$ and connects one side of $S_1$ to the other, implying that both copies of $S_1$ are in the boundary of the same connected component of $M\setminus\left(S_1\cup S_2\right)$.  Re-identifying the boundary copies of $S_1$, then, does not change the connectivity of $M\setminus\left(S_1\cup S_2\right)$ and yields the manifold $M\setminus S_2$.  But a single thrice-punctured sphere cannot separate an orientable three-manifold with torus boundary components, so $M\setminus S_2$ is connected. Since $M\setminus\left(S_1\cup S_2\right)$ and $M\setminus S_2$ have the same connectivity, this contradicts the fact that $\{S_1,S_2\}$ is a separating pair.  Thus a component of $\mathcal{A}$ that punctures the spheres $\{S_1,S_2\}$ must do so more than once. \end{proof}

Lemma \ref{lem:NotSP} will be used to show that longitudinal disks, the thrice-punctured sphere of Proposition \ref{prop:2PiDisk}, and reflection thrice-punctured spheres (except in the Borromean rings) do not contribute to separating pairs.  We begin with the case of longitudinal disks.

\begin{lemma}\label{lem:LongDisks}
If $C_1,C_2$ are distinct crossing circles of the FAL $\mathcal{A}$, then there is at most one longitudinal disk with punctures $C_1, C_2$.  In particular, a longitudinal disk is never part of a separating pair.
\end{lemma} 

\begin{proof}
Let $C_1,C_2$ be distinct crossing circles of the FAL $\mathcal{A}$, and let $p_1,p_2$ be the ideal vertices of $P_+$ projecting to cusps $C_1,C_2$.   Corollary \ref{cor:DisCDCircPackCrush} implies that, if a longitudinal disk $S_{\ell}$ has punctures $C_1,C_2$, then $p_1,p_2$ are two of the three points of tangency in a non-standard triple of the circle packing for $P_+$.  Given vertices $p_1, p_2$ of $P_+$, however, there is at most one non-standard triple containing them.  Indeed, if $p_1,p_2$ are vertices such a set, they both lie on one of the circles. The other two circles in the set must be those through $p_1, p_2$, so the vertices determine one such set.  Hence, there is at most one longitudinal disk with punctures $C_1,C_2$.

To see that separating pairs do not contain longitudinal disks, let $S_{\ell}$ be a longitudinal disk.  Since $S_{\ell}$ has punctures on distinct crossing circles, Lemma \ref{lem:NotSP} implies any thrice-punctured sphere that forms a separating pair with $S_{\ell}$ has the same punctures.  The previous paragraph, however, shows there is at most one such disk, so separating pairs do not contain longitudinal disks.
\end{proof}

\begin{figure}[h]
\[
\begin{array}{ccc}
\includegraphics{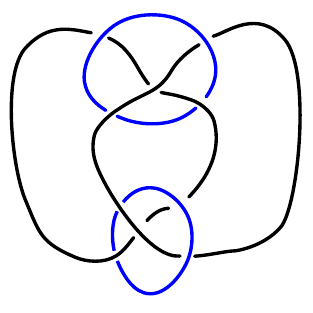}&\includegraphics{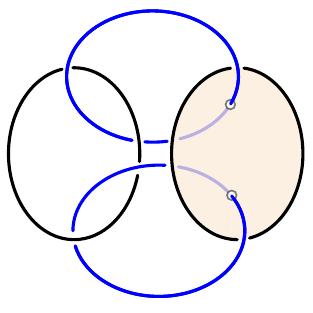}&\includegraphics{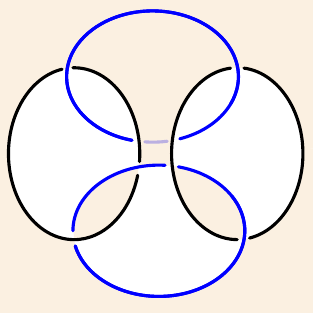}\\
(a)\textrm{ Fully Twisted}&(b)\textrm{ Obvious thrice-}&(b)\textrm{ Twisted-band}\\
\textrm{Borromean rings}&\textrm{punctured sphere}&\textrm{thrice-punctured sphere}
\end{array}
\]
\caption{The belted sum of $L_1$ and $L_2$}
\label{fig:Tet8}
\end{figure}

A second consequence of Lemma \ref{lem:NotSP} is that the fully twisted Borromean rings complement does not contain a separating pair.  Before proving this, note that the fully twisted Borromean rings is also the minimally twisted chain of four components (Figure \ref{fig:Tet8} parts $(a)$ and $(b)$ are isotopic).  Thrice-punctured spheres in this manifold (denoted $\mathbb{M}_4$ in \cite{Y2018}) are classified by Yoshida in \cite[Lemma 3.13]{Y2018}, and we recall Yoshida's classification here. Each subset of three components punctures two thrice-punctured spheres. Thus there are eight thrice-punctured spheres in total: the four ``obvious" twice-punctured disks of Figure \ref{fig:Tet8}$(b)$, together with those of Figure \ref{fig:Tet8}$(c)$ containing twisted bands.

\begin{lemma}\label{lem:Not2Pi3ps}
The fully twisted Borromean rings contains no separating pairs and is $b$-prime.  In particular, the thrice-punctured sphere with a single geodesic disk of area $2\pi$ is not part of a separating pair.
\end{lemma}

\begin{proof}
The fully twisted Borromean rings is the link depicted in Figure \ref{fig:Tet8} (see also \cite[Figure 6]{Y2018}).  Yoshida proves, \cite[Lemma 3.13]{Y2018}, its complement has exactly eight thrice-punctured spheres.  The spheres are the ``obvious" twice-punctured disks of Figure \ref{fig:Tet8}$(b)$, together with those of Figure \ref{fig:Tet8}$(c)$ containing twisted bands.  Note that each has boundary slopes along distinct cusps, and that each set of three distinct cusps contributes exactly two thrice-punctured spheres--those of Figure \ref{fig:Tet8}.  Two spheres with the same cusps intersect and do not form a separating pair.  Further, two thrice-punctured spheres that do not share the same cusps must have an isolated boundary slope and cannot separate by Lemma \ref{lem:NotSP}.  Thus the twisted Borromean rings do not contain a separating pair.
\end{proof}

We now show that reflection thrice-punctured spheres do not contribute to separating pairs in an FAL other than the Borromean rings.

\begin{lemma}\label{lem:Standard}
Let $\mathcal{A}$ be an FAL other than the Borromean rings, and let $M= S^3\setminus \mathcal{A}$ be its complement. No separating pair in $M$ contains a reflection thrice-punctured sphere. Moreover, two standard thrice-punctured spheres never form a separating pair for $M$.
\end{lemma}

\begin{proof}
We can assume $\mathcal{A}$ is not the fully twisted Borromean rings, since Lemma \ref{lem:Not2Pi3ps} implies the result in this case. We argue that a reflection thrice-punctured sphere cannot form a separating pair with a non-reflection disk.  If $S$ is a reflection thrice-punctured sphere, the punctures of $S$ along knot circles are (generalized) longitudes.  Additionally, let $m, \ell$ denote the meridian and longitude slopes of a crossing circle $C$ that punctures $S$. If $C$ is a flat crossing circle it punctures $S$ in two meridians, and if $C$ is twisted it punctures $S$ in the $\ell + 2m$ slope.  Theorem \ref{thm:NR3ps} implies that punctures of a non-reflection disk are either meridians of knot circles, longitudes of crossing circles, or meridians of crossing circles, all of which intersect the punctures of $S$. Thus any non-reflection disk that shares a puncture with $S$ intersects it, and cannot form a separating pair with $S$.

Now assume $\{S,T\}$ are both reflection thrice-punctured spheres that form a separating pair.  Since $M\setminus (S\cup T)$ is disconnected, and no proper subset of the reflection surface separates $M$, the reflection surface must be $S\cup T$.  We show that the reflection surface equaling $S\cup T$ implies that $\mathcal{A}$ is the Borromean rings, whose reflection surface does form a separating pair. Now the nerve of the circle packing for $\mathcal{A}$ is a triangulation of $\mathbb{S}^2$, so it has at least four vertices.  Vertices of the nerve correspond to unshaded faces so $P_+$ has at least four unshaded faces, each of which is an ideal polygon.  Since they glue up to two thrice-punctured spheres, area considerations imply the unshaded faces of $P_+$ must be four ideal triangles.   Therefore, $P_+$ is a regular ideal octahedron and the only FALs resulting from gluing two regular ideal octahedra are the Borromean rings, possibly with twisted crossing circles.  Twisting one crossing circle in the Borromean rings makes one of the reflection surface components non-orientable, so it cannot be a thrice-punctured sphere, and we've already assumed $\mathcal{A}$ is not the fully twisted Borromean rings.  Thus the only FAL with two thrice-punctured spheres for a reflection surface is the Borromean rings themselves.

To finish the proof it remains to show that two standard crossing disks do not form a separating pair.  However, distinct standard crossing disks are punctured by distinct crossing circle components, and Lemma \ref{lem:NotSP} shows they do not form a separating pair.
\end{proof} 

The following theorem characterizes separating pairs in FAL complements, summarizing the results of this section.

\begin{thm}\label{thm:GeneralBS}
Let $\mathcal{A}$ be an FAL other than the Borromean rings with complement $M= S^3\setminus \mathcal{A}$. Suppose that $S_1,S_2$ are disjoint thrice-punctured spheres in $M$. 

The pair $\{S_1,S_2\}$ is a separating pair if and only if each is either a crossing disk or a singly-separated disk, and their longitudinal slopes coincide.
\end{thm}

\begin{proof}
Suppose $\{S_1,S_2\}$ is a separating pair.  Proposition \ref{prop:2PiDisk} and Lemma \ref{lem:Not2Pi3ps} imply that neither $S_1$ nor $S_2$ has an area $2\pi$ geodesic disk. Lemmas \ref{lem:LongDisks} and \ref{lem:Standard} show that neither $S_1$ nor $S_2$ are longitudinal disks or reflection thrice-punctured spheres.   After eliminating these possibilities, crossing disks (standard or non-standard) and singly-separated disks are the only remaining types for $S_1$ and $S_2$.  Now both crossing and singly-separated disks have a unique longitudinal slope (which is along a crossing circle).  If the longitudinal slopes of $S_1$ and $S_2$ are along distinct crossing circle cusps of $M$, then Lemma \ref{lem:NotSP} implies the pair $\{S_1,S_2\}$ is not a separating pair.  Thus the longitudinal slopes of $S_1$ and $S_2$ must coincide.

Conversely, let $S_1$ and $S_2$ be two crossing or singly-separated disks that share a longitudinal puncture along the same crossing circle $C$.  One obtains an embedded two sphere $\mathbb{S}^2$ in $\mathbb{S}^3$ from $S_1\cup S_2$ by including the crossing circle $C$, together with the points of $\mathcal{A}$ corresponding to meridional punctures of $S_1\cup S_2$.  Since $\mathbb{S}^2$ separates $\mathbb{S}^3$, $M\setminus \left(S_1\cup S_2\right)$ is also disconnected, and $\{S_1,S_2\}$ is a separating pair.
\end{proof}

As observed earlier, Theorem \ref{thm:GeneralBS} shows that every separating pair in an FAL complement consists of a shared crossing-circle longitudinal puncture together with meridional punctures.  Thus every punctured sphere decomposition in an FAL complement can be realized as a belted-sum decomposition. 

Separating pairs consisting of two crossing disks will play an important roll in what follows.  Corollary \ref{cor:PrimeIs1Disk}, for example, shows that an FAL is $b$-prime if and only if there are no separating pairs with two crossing disks.  Moreover, such separating pairs are exactly those used to define canonical belted sum decompositions of FALs (see Definition \ref{defn:CanonicalBS}).  For these reasons we now prove that every separating pair consisting of two crossing disks decomposes an FAL into the belted sum of two FALs.

\begin{lemma}\label{lem:CrossDiskSepPair}
Let $\{S_1,S_2\}$ be a separating pair consisting of two crossing disks in the FAL complement $M=\mathbb{S}^3\setminus \mathcal{A}$, and let $M = M_1\#_b M_2$ denote a belted-sum decomposition they determine.  Then both $M_1$ and $M_2$ are FAL complements. 
\end{lemma}

\begin{proof}
Let $\{S_1,S_2\}$ be a separating pair consisting of two crossing disks and $M = M_1\#_b M_2$ their corresponding belted-sum decomposition.  Consider the case where neither $S_i$ is standard (the case where one of the $S_i$ is the standard crossing disk is similar). Then $S_1$ and $S_2$ share a longitudinal puncture along the same crossing circle $C$ of $\mathcal{A}$.  

We begin by slicing a fundamental polyhedron for $M$ into two polyhedra and determining face pairings in the individual pieces. Then show these face pairings result in FALs.  This is the reverse of Adams proof of \cite[Theorem 4.5]{ad1}, where he glues fundamental polyhedra for the individual manifolds along faces corresponding to chosen thrice-punctured spheres.  The standard polyhedral decomposition for FALs facilitates our process for us and, since $P_-$ is the reflection of $P_+$, our description focuses on $P_+$.

The crossing circle $C$ shared by $\{S_1,S_2\}$ corresponds to exactly one ideal vertex of $P_+$ and we assume, possibly after a M\"obius transformation, that it corresponds to infinity (using the upper-half space model for $\mathbb{H}^3$).  Figure \ref{fig:CrossDiskFALBS}$(a)$ illustrates $P_+$ with geodesic disks $S_i'$ for the non-standard crossing disks $S_i$.  When viewed from above, $P_+$ is rectangular and the vertical shaded faces (labeled $T'$) project to the standard crossing disk for $C$.  The geodesic disks $S_i'$ are ideal triangles parallel to the $T'$ since all have longitudinal punctures along $C$. The shaded interstices of Figure \ref{fig:CrossDiskFALBS}$(a)$ may or may not contain additional circles of the circle packing. 

\begin{figure}[h]
\[
\begin{array}{ccc}
\includegraphics{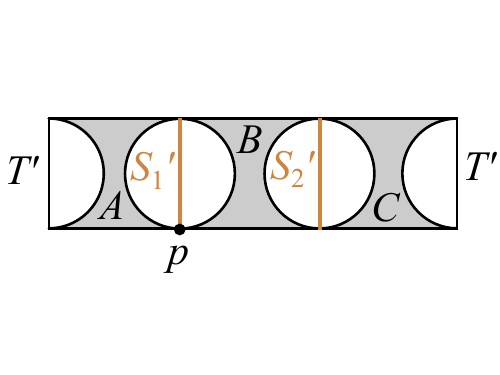}&\hspace{0.2in}&\includegraphics{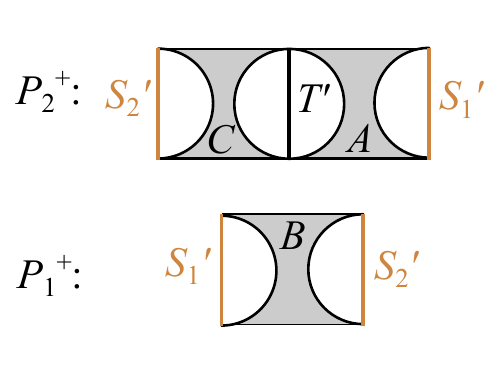}\\
(a)\textrm{ Separating pair in $P_+$}&&(b)\textrm{ The polyhedra $P_i^+$}
\end{array}
\]
\caption{Belted-sum decomposition if both $\{S_1,S_2\}$ are non-standard crossing disks}
\label{fig:CrossDiskFALBS}
\end{figure}

Now slice $P_+$ along the $S_i'$ to obtain three pieces, then glue the outside pieces along the shaded triangles $T'$.  The result is two polyhedra, and let $P_1^+$ denote the one between the $S_i'$ while $P_2^+$ denotes that formed from the outsides (see Figure \ref{fig:CrossDiskFALBS}$(b)$).  A similar slicing along $P_-$ yields the reflections $P_1^-$ and $P_2^-$.  Note that in the case where $T'$ projects to a twisted crossing disk, the polyhedron $P_2^+$ will consist one outside piece from each of $P_{\pm}$.  

Now consider the gluing maps on faces of $P_i^{\pm}$. The gluing maps of $P_{\pm}$ identify unshaded faces of $P_i^{\pm}$ to each other since the disks $S_i'$ are orthogonal to the reflection surface. Moreover, since the $S_i$ are crossing disks, exactly one vertex of the ideal triangles $S_i'$ corresponds to a crossing circle.  This means the triangular faces $T'$ are the only ones of $P_{\pm}$ adjacent to $S_i'$ that are identified on $P_{\pm}$, and they are glued to form $P_2^{\pm}$. The face pairings on the remaining shaded triangles of $\partial P_{\pm}$ do not connect $P_1^{\pm}$ to $P_2^{\pm}$. Thus all faces of $P_i^{\pm}$ inherit gluing maps from $P_{\pm}$ except the shaded triangles corresponding to the $S_i'$, and those maps do not connect the $P_1^{\pm}$ to  $P_2^{\pm}$. Pairing the copies of $S_i'$ of $P_i^+$ along their shared vertex completes the gluing instructions so that unshaded faces glue to their reflections and shaded faces glue to their mate or its reflection.   

Indentifying the $S_i'$ either way yields a so-called \emph{admissible} gluing pattern from \cite{hot}, and \cite[Theorem 3.2 ]{hot} guarantees that $P_i^{\pm}$ with these face-pairings yield complete hyperbolic manifolds $M_i$. One can then form the circle packing and nerve corresponding to these polyhedral decompositions, and apply \cite[Lemma 2.4]{pu1} to conclude the manifolds are FAL complements.  Let $M_i$ denote the FAL resulting from gluing $P_i^{\pm}$.   The only ambiguity in this process is whether the $S_i'$ project to flat or twisted crossing disks. 

Beginning with $M_1$ and $M_2$ and forming the belted sum $M_1 \#_b M_2$ along the disks $D_i$, as in \cite{ad1}, reverses this process on the polyhedral level, thereby yielding $M$ by construction.  

\begin{figure}[h]
\begin{center}
\includegraphics{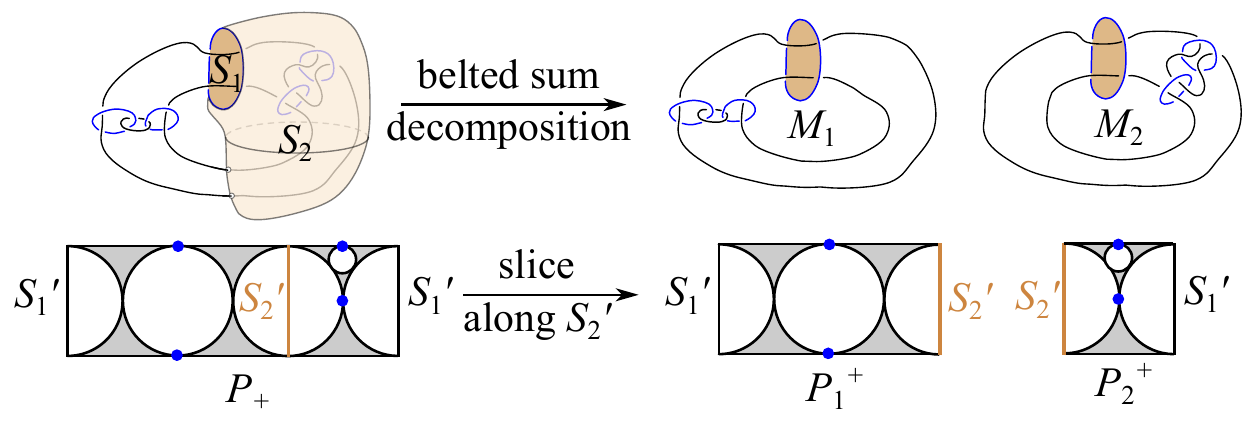}
\end{center}
\caption{$\{S_1,S_2\}$ with one standard and one non-standard crossing disk}
\label{fig:PolyBS}
\end{figure}

The case where the disk $S_1$ of the pair $\{S_1,S_2\}$ is a standard crossing disk is similar.  The difference is that one only slices $P_+$ along $S_2'$ (respectively $P_-$ along the reflection of $S_2'$), yielding only two pieces which are the $P_1^{\pm}$ and $P_2^{\pm}$ described (see Figure \ref{fig:PolyBS}).  
\end{proof}

\begin{figure}[h]
\begin{center}
\includegraphics{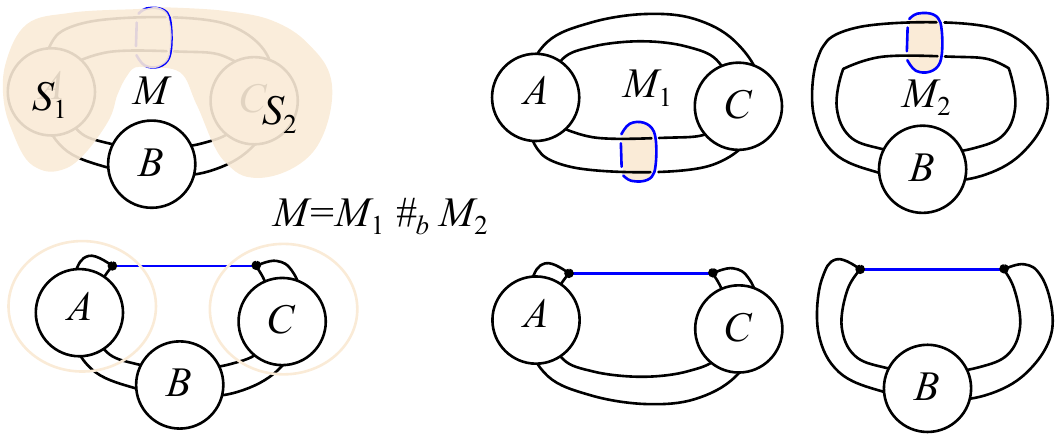}
\end{center}
\caption{FAL belted sum and crushtaceans}
\label{fig:BSTwoCD}
\end{figure}

The case just argued, where both $S_i$ are non-standard crossing disks, can be illustrated in the crushtacean as well (see Figure \ref{fig:BSTwoCD}). The tangles $A,B,C$ of Figure \ref{fig:BSTwoCD} correspond to the interstices of the same label in Figure \ref{fig:CrossDiskFALBS}.

Thus every separating pair consisting of two crossing disks decomposes an FAL complement into the belted sum of two simpler FAL complements.  We turn our attention to the questions of $b$-primality and canonical decompositions.

\section{On $b$-Prime FALs}\label{sec:FALPrime} 

This section provides three characterizations of $b$-prime FALs. The Borromean rings and its half-twist partners are a special case, and considered separately. For FALs with more than three crossing circles, Theorem \ref{thm:CrushCharPrime} shows that an FAL is $b$-prime if and only if each non-trivial three-edge cut in its crushtacean is fully painted.  Interpreting this in terms of crossing disks shows that an FAL is $b$-prime if and only if each crossing circle bounds a unique crossing disk (see Corollary \ref{cor:PrimeIs1Disk}).  These results follow readily from those of Section \ref{sec:SP}. 
After introducing diagrammatic terminology, Theorem \ref{thm:PrimeIsTrivialOrbit} shows that $b$-prime FALs come from fully augmenting diagrams with trivial flype orbits. 

First consider the special case of FALs with two crossing circles--the Borromean rings and its half-twist partners.   

\begin{lemma}\label{lem:NonPrimeBrings}
The fully twisted Borromean rings complement is $b$-prime. The complement of the Borromean rings with both, or exactly one, flat crossing circle are $b$-composite.  

In fact, both are the belted sum of two Whitehead link complements and both have a separating pair consisting of one crossing and one singly-separated disk.
\end{lemma}

\begin{proof}
We've already seen in Lemma \ref{lem:Not2Pi3ps} that the fully twisted Borromean rings is $b$-prime.  Moreover, Adams observed (see \cite{ad1}) that the Borromean rings complement is the belted sum of two copies of the Whitehead link complement.  In our terminology, the singly-separated disk of Figure \ref{fig:BRings2WL} forms a separating pair with the top crossing disk. This pair, with the appropriate choice of gluing homeomorphisms, decomposes the Borromean rings into two copies of the Whitehead link.

\begin{figure}[h]
\begin{center}
\includegraphics{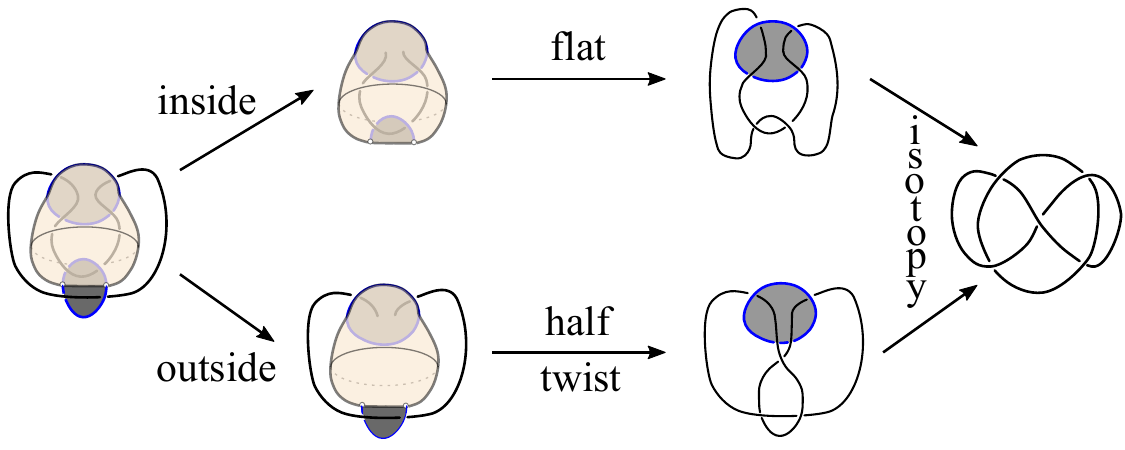}
\end{center}
\caption{Borromean rings as belted sum of two Whitehead links}
\label{fig:BRings2WL}
\end{figure}

Finally, the separating pair of the original Borromean rings is preserved when adding a half-twist to the top crossing disk; therefore, the complement of the Borromean rings with one twisted crossing circle is $b$-composite.  As in the case of the Borromean rings complement, it is the belted sum of two Whitehead link complements, with appropriate choice of gluing homeomorphisms. 
\end{proof}

We now treat the general case, characterizing $b$-prime FALs in terms of their painted crushtaceans.

\begin{thm}\label{thm:CrushCharPrime}
Let $\mathcal{A}$ be a fully augmented link with at least three crossing circles, and let $\Gamma^{\ast}$ be its crushtacean.  

The complement $M = \mathbb{S}^3\setminus\mathcal{A}$ is $b$-prime if and only if every three-edge cut of $\Gamma^{\ast}$ is fully painted.
\end{thm}

\begin{proof}
First recall that every non-trivial three-edge cut has either one or three painted edges (see discussion near Figure \ref{fig:NSMT}).  If $\Gamma^{\ast}$ has a once-painted, non-trivial three-edge cut, then $M$ has a non-standard crossing disk by Corollary \ref{cor:DisCDCircPackCrush}.  This forms a separating pair with the standard disk sharing the crossing circle puncture, by Theorem \ref{thm:GeneralBS}, and $M$ is not $b$-prime.

Conversely, suppose every three-edge cut of $\Gamma^{\ast}$ is fully painted, so that every non-standard $n$-disk is longitudinal.  Theorem \ref{thm:GeneralBS}, then, implies that every separating pair must contain a singly-separated disk.  We now argue that if $M$ contains a singly-separated disk, it contains a non-standard crossing disk, contradicting the fact that $n$-disks in $M$ are assumed longitudinal.

\begin{figure}[h]
\[
\begin{array}{ccc}
\includegraphics{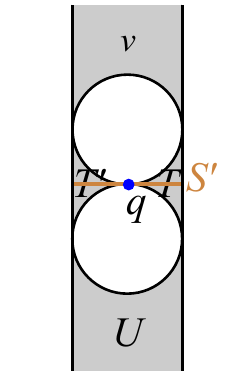}&\hspace{0.25in}&\includegraphics{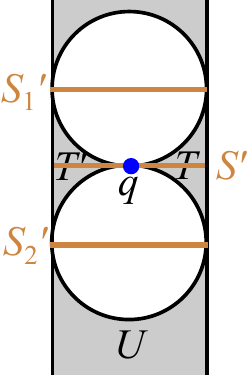}\\
(a)  \textrm{ $S'$ projects to $S$}&&(b) \textrm{ $S_i'$ neighbors of $S'$}
\end{array}
\]
\caption{Ideal triangular neighbors of material-vertex disks}
\label{fig:NeighborsOfMVD}
\end{figure}

Let $S$ be a singly-separated disk in an FAL complement $M$, and $S'$ one of its material vertex disks. Proposition \ref{prop:3psSeparatesD} shows that circle packing near $S'$ is as illustrated in Figure \ref{fig:NeighborsOfMVD}$(a)$.   If both interstices $U$ and $v$ are trivial, this is the circle packing for an FAL with two crossing circles. Since $M$ has at least three crossing circles, at least one of the interstices, say $U$, has additional circles in the circle packing.  Now $S'$ has two neighboring ideal triangles labeled $S_i'$ in Figure \ref{fig:NeighborsOfMVD}$(b)$, and we focus on $S_2'$.  Since $q$ is the only vertex of the triangles $T, T'$ corresponding to a crossing circle, the ideal vertices at the base of $S_2'$ correspond to knot circles in $\mathcal{A}$.  Proposition \ref{prop:3psDisjointFromCD} implies that $S_2'$ projects to a crossing disk and, since interstice $U$ is non-trivial, the disk $S_2'$ is not a face of $P_+$.  Thus $S_2'$ projects to a non-standard crossing disk, contradicting the assumption that every three-edge cut is fully painted.

If every three-edge cut of $\Gamma^{\ast}$ is fully painted, then, $M$ contains only standard crossing disks and is void of singly-separated disks.  Theorem \ref{thm:GeneralBS} implies $M$ is void of separating pairs, and is $b$-prime.
\end{proof}

The correspondence between once-painted, three-edge cuts of the crushtacean and non-standard crossing disks immediately leads to the following characterization of $b$-prime FALs.

\begin{cor}\label{cor:PrimeIs1Disk}
Let $\mathcal{A}$ be a fully augmented link with at least three crossing circles.

The complement $M = \mathbb{S}^3\setminus\mathcal{A}$ is $b$-prime if and only if every crossing circle of $\mathcal{A}$ bounds a unique crossing disk.
\end{cor}

\begin{proof}
Every crossing circle of $\mathcal{A}$ bounds a standard crossing disk. Corollary \ref{cor:DisCDCircPackCrush} implies that non-standard crossing disks are in one-to-one correspondence with once-painted three-edge cuts of the crushtacean.  Thus every crossing circle of $\mathcal{A}$ bounds a unique crossing disk (i.e. the standard one) if and only if there are no once-painted three-edge cuts of the crushtacean. The result follows from Theorem \ref{thm:CrushCharPrime}.
\end{proof}

Observe that the criteria of Corollary \ref{cor:PrimeIs1Disk} can be rephrased to say that $M = \mathbb{S}^3\setminus\mathcal{A}$ contains no non-standard crossing disks.

We now interpret FAL $b$-primality in terms of the link diagram that generated the FAL.  The main result is that an FAL $\mathcal{A}_D$ obtained by augmenting a diagram $D$ is $b$-prime if and only if it generates only trivial flype orbits.  We review the necessary diagrammatic properties involved before proving the theorem, and refer the reader to \cite{ca} for further diagrammatic details.

\begin{figure}[h]
\begin{center}
\includegraphics{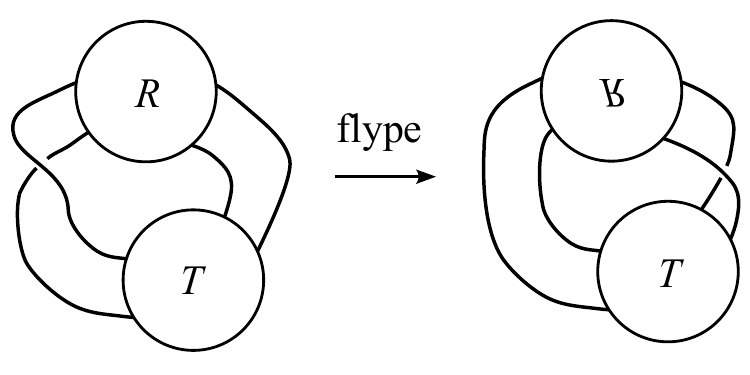}
\end{center}
\caption{ Flype on a link diagram}
\label{fig:Flype}
\end{figure}

An \emph{m-tangle} in a regular link diagram $D$ is a disk $\mathcal{D}$ in the projection plane whose boundary intersects $D$ in $m$ points.  An $m$-tangle is trivial if it contains no crossings of $D$.  A flype changes one diagram of a link $L$ into another by rotating a 4-tangle $180^{\circ}$, resulting in a crossing moving across the tangle as in Figure \ref{fig:Flype}.  A flype is non-trivial if neither of the 4-tangles involved are trivial. Flypes were introduced by Tait during his original efforts to enumerate links in the late 1800's, when he conjectured that any two reduced alternating diagrams of a link are related by a sequence of flypes.  This conjecture was finally proven by Menasco and Thistlethwaite in \cite{mt}.

A \emph{flype tangle} for a crossing is a 4-tangle that cannot be decomposed into smaller flype-admitting 4-tangles for that crossing.  A flype orbit for a crossing is a collection of flype tangles for the crossing, and a flype orbit is trivial if it consists of a single flype tangle.  Each crossing generates two flype orbits, and Calvo shows that at most one of them is non-trivial (see \cite[Lemma 4]{ca}).

\begin{center}
\begin{figure}[h]
\[
\begin{array}{ccccc}
 \includegraphics[width=1.25in]{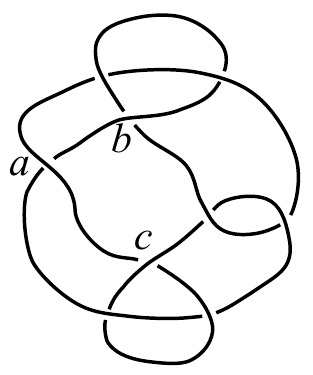}&\hspace{0.1in}&\includegraphics[width=1.25in]{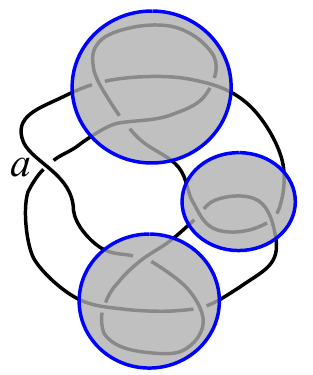}&\hspace{0.1in}&\includegraphics[width=1.25in]{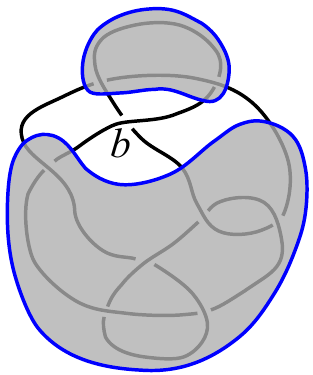} \\
(a)\textrm{ Knot }9_{28}&& (b)\ a\textrm{ flype orbit } && (c)\ b\textrm{ flype orbit }
\end{array}
\]
\caption{A knot and its flype orbits}
\label{fig:FlypeOrbit}
\end{figure}
\end{center}

Figure \ref{fig:FlypeOrbit} illustrates these definitions.  Knot $9_{28}$ of Rolfsen's tables has six twist regions (see Figure \ref{fig:FlypeOrbit}$(a)$).  The twists with nontrivial flype orbits consist of single crossings, while the 2-crossing twists have only trivial flype orbits.  Parts $(b)$ and $(c)$ of Figure \ref{fig:FlypeOrbit} depict flype tangles and orbits for crossings $a$ and $b$, respectively (the flype orbit for crossing $b$ is more easily seen after an isotopy of the projection sphere).  By symmetry, rotating the diagram of Figure \ref{fig:FlypeOrbit}$(c)$ about a horizontal axis yields the flype tangles and orbit of crossing $c$.

Let $D$ be a reduced alternating diagram for a link $L$, and recall that fully augmenting $L$ from $D$ amounts to choosing a crossing circle around each twist in $D$. There is a unique choice of crossing circle associated with each twist of $D$ containing more than one crossing, and the corresponding flype orbit may or may not be trivial.  Single-crossing twists in $D$ admit two choices for crossing circles, and at most one of the choices corresponds to a non-trivial flype orbit by \cite[Lemma 4]{ca}. Given an FAL $\mathcal{A}_L$ arising from $D$, the \emph{flype orbits induced by} $\mathcal{A}_L$ are those determined by the crossing circles. 

\begin{center}
\begin{figure}[h]
\[
\begin{array}{ccc}
 \includegraphics[width=1.1in]{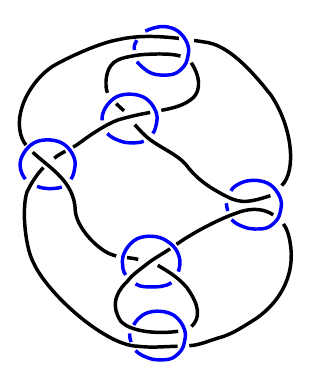}\includegraphics[width=1.1in]{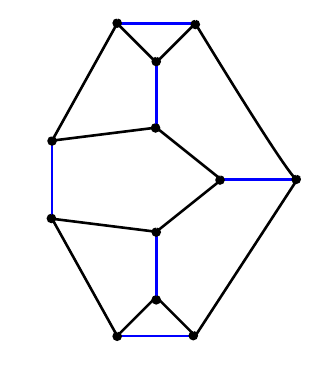}&&\includegraphics[width=1.1in]{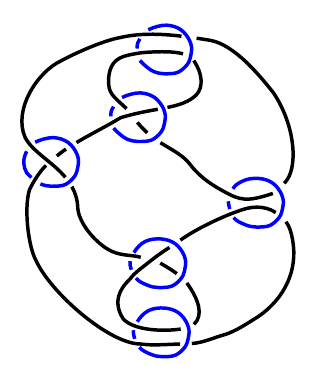}\includegraphics[width=1.1in]{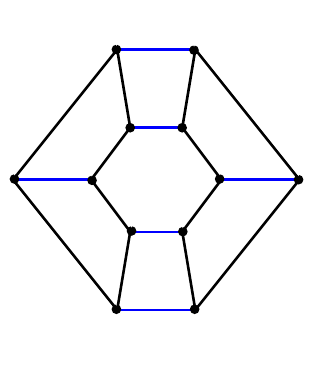}\\
(a)\ \mathcal{A}\textrm{ and its crushtacean}&& (b)\ \mathcal{A}'\textrm{ and its crushtacean}
\end{array}
\]
\caption{$b$-composite and $b$-prime augmentations of knot $9_{28}$}
\label{fig:Aug9_28}
\end{figure}
\end{center}

For example, the augmentation $\mathcal{A}'$ of Figure \ref{fig:Aug9_28}$(b)$ is obtained by changing the choice of crossing circle for every single-crossing twist in Figure \ref{fig:Aug9_28}$(a)$.  The FAL $\mathcal{A}'$ induces all trivial flype orbits to augment $D$, while $\mathcal{A}$ induces non-trivial flype orbits for all single-crossing twists. The crushtaceans of Figure \ref{fig:Aug9_28} illustrate the fact that the same diagram can have a $b$-prime augmentation as well as $b$-composite ones. Note that every painted edge corresponding to a twisted crossing circle in Figure \ref{fig:Aug9_28}$(a)$ is contained in at least one non-trivial, once-painted three-edge cut; whereas the crushtacean for $\mathcal{A}'$ contains no non-trivial, three-edge cuts.

We now have the terminology needed to give the following diagrammatic characterization of $b$-prime FALs.  A technical note in the following theorem: the assumption that $D$ has more than two twists prevents $\mathcal{A}$ from being the Borromean rings or one of its half-twist partners.

\begin{thm}\label{thm:PrimeIsTrivialOrbit}
Let $D$ be a twist-reduced, non-split, prime link diagram with more than two twist regions, and let $\mathcal{A}$ be a full augmentation of $D$.  

The fully augmented link $\mathcal{A}$ is $b$-prime if and only if it induces only trivial flype orbits.
\end{thm}

\begin{proof}
The assumptions on the diagram $D$ ensure that $\mathcal{A}$ is a hyperbolic FAL other than the Borromean rings or its half-twist partners (See \cite[Theorem 2.5]{pu1}).  The key observation to prove is that the manifold $M=\mathbb{S}^3\setminus \mathcal{A}$ contains non-standard crossing disks if and only if $\mathcal{A}$ induces a non-trivial flype orbit. The fact that the changes made in passing from the diagram $D$ to the FAL $\mathcal{A}$ are local allows us to argue diagramaticly. The result then follows from Corollary \ref{cor:PrimeIs1Disk}.

Suppose the augmentation $\mathcal{A}$ of $D$ induces a non-trivial flype orbit for a twist region $\tau$ with the chosen crossing circle $C$.  Let $\mathcal{T}_1,\dots,\mathcal{T}_n$, with $n \ge 2$, denote the flype tangles in the order encountered traversing the orbit for $\tau$. One can choose pairwise disjoint three-ball neighborhoods $\mathbb{B}_i$ of the $\mathcal{T}_i$ such that $\cup_{i=1}^n \mathbb{B}_i$ contains all components of $\mathcal{A}$ except for $C$ and the strands of $D$ connecting the $\mathcal{T}_i$. In addition, choose the $\mathbb{B}_i$ so that they are disjoint from the crossing circle $C$ and its standard crossing disk $S$.  The crossing circle $C$ bounds a disk $S_i$ disjoint from $\cup_{i=1}^n \mathbb{B}_i$ that intersects the arcs of $D$ joining $\mathcal{T}_i$ and $\mathcal{T}_{i+1}$.  Each $S_i$ is a thrice-punctured sphere with a longitudinal puncture along $C$ and two meridional knot circle punctures, so $S_i$ is a non-standard crossing disk.  Moreover, the $S_i$ are distinct (not isotopic) since the $\mathbb{B}_i$ contain components of $\mathcal{A}$ which separate them.  Hence if $\mathcal{A}$ induces a non-trivial flype orbit, $M$ contains non-standard crossing disks and is $b$-composite.

Conversely, suppose $C$ is a crossing circle of $\mathcal{A}$ that bounds a non-standard crossing disk $S'$ in $M$, and let $S$ denote the standard crossing disk for $C$.  Let $\tau$ be the twist region of the diagram $D$ augmented by $C$, we wish to show that $\tau$ has a non-trivial flype orbit.  By Theorem \ref{thm:GeneralBS}, the pair $\{S,S'\}$ is a separating pair.  The non-separating geodesics of $\{S,S'\}$ are contained in the projection plane of $D$ and form, together with points of their punctures, a simple closed loop $\gamma$ in $D$.  The projection plane intersects the pair $\{S,S'\}$ in their non-separating geodesics even when $C$ is a twisted crossing circle.  This follows from the fact that the half-twist changes the reflection surface and not the projection plane containing $D$.

Now $\gamma$ is the simple closed curve in the projection plane determined by non-separating geodesics of $\{S,S'\}$ together with the necessary puncture points.  The boundary curves of an appropriately chosen product neighborhood $\gamma \times (-\epsilon,\epsilon)$, then, form 4-tangles containing non-trivial portions of $D$, since $\{S,S'\}$ is a separating pair.  These tangles allow for a non-trivial flype, showing the flype orbit of $\tau$ is non-trivial.

Since $\mathcal{A}$ induces non-trivial flype orbits if and only if $M$ contains non-standard crossing disks, the proof is complete by Corollary \ref{cor:PrimeIs1Disk}.
\end{proof}

We now have three characterizations of $b$-prime FALs: one combinatorial, another geometric, and the third diagrammatic.  The next section applies this to create canonical belted-sum decompositions of an FAL into $b$-prime FALs together with copies of the Borromean rings.

\section{Canonical FAL Decompositions}\label{sec:CanonDec}

There are two desirable properties for belted-sum decompositions of an FAL complement: first, summands ought to be $b$-prime and, second, summands should be FALs to maintain the tractable geometry of FALs in the process.  These two properties can nearly be achieved simultaneously.  More precisely, this section introduces a canonical method for decomposing an FAL in which every summand is an FAL complement and, if a summand is not $b$-prime, then it is a copy of the Borromean rings with at least one flat crossing circle.  If desired, each summand can be made $b$-prime by decomposing each copy of the Borromean rings into two Whitehead links by Lemma \ref{lem:NonPrimeBrings}. The result can then be recast as follows: there is a canonical belted-sum decomposition of an FAL complement in which each summand is either a $b$-prime FAL or a copy of the Whitehead link.  The decomposition is canonical in the sense that the prescribed choices in the process produce a unique decomposition.

A technical lemma is required before defining the canonical belted sum decomposition, which is that distinct crossing disks in an FAL complement are disjoint.  This follows from the more general result that $n$-disks are disjoint in an FAL complement.

\begin{lemma}\label{lem:CrossDiskDisjoint}
If $M = \mathbb{S}^3\setminus\mathcal{A}$ is an FAL complement, then distinct $n$-disks in $M$ are disjoint. 
\end{lemma}
\begin{proof}
Any two $n$-disks $S_i$ and $S_j$ are both orthogonal to the reflection surface $R$ of $M$, and both intersect $R$ in their non-separating geodesics (Proposition \ref{prop:3psDisjointFromCD}).  If $S_i\cap S_j$ contains a non-separating geodesic $\gamma$, then $S_i=S_j$ since they both are orthogonal to $R$ along $\gamma$. Moreover, distinct $n$-disks cannot intersect in a geodesic that is separating on one and non-separating on the other, since separating geodesics on $n$-disks are not contained in $R$ while non-separating geodesics are.  Thus $S_i\cap S_j$ can only contain a geodesic that is separating on both. Yoshida, however, has classified intersections of two thrice-punctured spheres in an orientable hyperbolic three manifolds (see \cite[Proposition 3.1]{Y2018}), showing that they cannot intersect in a geodesic that separates both.  Hence $n$-disks are pairwise disjoint.
\end{proof}

In particular, Lemma \ref{lem:CrossDiskDisjoint} shows that crossing disks in an FAL complement $M$ are disjoint, and we can slice $M$ along any or all of them simultaneously.  Recall that two crossing disks $\{S_1,S_2\}$ form a separating pair if and only if they share the same longitudinal crossing circle puncture.  Also, the proof of Lemma \ref{lem:CrossDiskSepPair} shows that such a pair decomposes $M$ into two FAL summands, and the only ambiguity is whether the corresponding crossing circles in the summands are twisted or flat. In defining the canonical decomposition, we make the convention to identify copies of the $S_i$ within $P_+$, resulting in flat crossing circles. We are now prepared to define the canonical belted sum decomposition of an FAL.

\begin{defn}\label{defn:CanonicalBS}
Let $\mathcal{A}$ be an FAL with complement $M=\mathbb{S}^3\setminus\mathcal{A}$.  The \emph{\textbf{canonical belted sum decomposition}} of $M$ is obtained by the following procedure:

\begin{enumerate}
\item First slice $M$ along all crossing disks that belong to a separating pair, then
\item Identify boundary spheres in the resulting components which correspond to separating pairs of $M$ so as to create flat crossing disks in the summands. 
\end{enumerate}
\end{defn}

Two component FALs are their own canonical decomposition, since each crossing circle bounds a unique crossing disk.
More generally, Corollary \ref{cor:PrimeIs1Disk} shows that if $M$ is already $b$-prime then each crossing circle bounds a unique crossing disk, so the canonical belted sum decomposition of a $b$-prime FAL is itself. 

The canonical belted sum decomposition is well-defined.  Indeed, since all crossing disks are disjoint by Lemma \ref{lem:CrossDiskDisjoint}, the result of Step $(1)$ is a collection of manifolds with thrice-punctured spheres in their boundaries. Moreover, boundary spheres are paired along their common longitudinal crossing circle punctures, uniquely determining which boundary spheres are to be identified.  Finally, the restriction that the crossing circles resulting from these identifications be flat uniquely determines the summands of the canonical decomposition. Thus there is a unique canonical decomposition for each given FAL.

We now prove that summands in canonical decompositions are either $b$-prime FALs or copies of the Borromean rings with zero or one twisted crossing circles. To prove this, we must show that decomposing an FAL complement along crossing disks can be done simultaneously, that the summands are always FALs, and that each summand that is not is $b$-prime is a copy of the Borromean rings with zero or one twisted crossing disk(s).

\begin{thm}\label{thm:CanonIsBPrime}
Let $\mathcal{A}$ be an FAL with at least three crossing circles and let $M=\mathbb{S}^3\setminus\mathcal{A}$ be its complement. 

Each summand in the canonical belted sum decomposition of $M$ is either a $b$-prime FAL complement or the complement of the Borromean rings with at least one flat crossing circle.
\end{thm}

\begin{proof}
Let $M = M_1\#_b \dots \#_b M_m$ be the canonical belted sum decomposition of $M$ and let $S_1,\dots,S_n$ denote all crossing disks (both standard and non-standard) involved in a separating pair of $M$. As noted above, the slicing operation $M\setminus \left(\cup_{i=1}^nS_i \right)$ can be done simultaneously since distinct crossing disks in $M$ are disjoint. 

Now suppose $C$ bounds the crossing disks $S_{1}^C,\dots, S_{k}^C$, with $k\ge 2$.  Further, suppose they are labeled in the order encountered along a meridian of $C$, with $S_{1}^C$ denoting the standard crossing disk.  Any two of them form a separating pair, so the two copies of $S_{i}^C$ are on distinct components of $M\setminus \left(\cup_{i=1}^nS_i \right)$. In addition, (at least) one component of $M\setminus \left(\cup_{i=1}^nS_i \right)$ has one copy of each consecutive $S_{i}^C,S_{{i+1}}^C$ in its boundary, and they share a longitudinal puncture along $C$.  

Let $M_j'$ denote the component of $M\setminus \left(\cup_{i=1}^nS_i \right)$ that forms the summand $M_j$ when boundary spheres are identified. The previous paragraph implies the boundary of $M_j'$ consists of an even number of thrice-punctured spheres which are paired according to their longitudinal punctures in the original manifold $M$. Preserving longitudinal and meridional punctures while identifying paired thrice-punctured spheres in $M_j'$ yields a FAL $M_i$ by Lemma \ref{lem:CrossDiskSepPair}. Our initial goal is to show that each $M_j$ with at least three crossing circles is $b$-prime.

Our first task is to show that crossing disks of $M$ contribute a unique crossing disk for each crossing circle of $M_j$, which we choose to be standard.  First note that if $C$ is crossing circle in $M$ that bounds a unique crossing disk $S_C$, then the belted sum operations in the canonical decomposition of $M$ are disjoint from a neighborhood of $C\cup S_C$.  Thus $C$ bounds the disk $S_C$ in $M_j$.  Now suppose that $C$ bounds at least one non-standard crossing disk in $M$.  Then copies of consecutive crossing disks in $M$ along $C$ are glued to form a disk $S_C$ for $C$ in $M_j$.  In either case, the crossing disks of $M$ contribute a unique crossing disk $S_C$ for each $C$ in $M_j$, and we choose these to be the standard disks in $M_j$.  

We now argue that the process of canonical decompositions do not create additional crossing disks in the summands.  Lemma \ref{lem:CrossDiskSepPair} implies that each summand $M_j$ is an FAL, as each can be formed by consecutive belted sums along separating pairs of crossing disks. If $M_j$ is $b$-composite, then Corollary \ref{cor:PrimeIs1Disk} implies there is a crossing circle $C$ of $M_j$ that bounds a non-standard crossing disk $S$.  Note that $S$ is not part of a belted sum operation needed to build $M$, since only chosen standard crossing disks of $M_j$ are used to reconstruct $M$.  By Lemma \ref{lem:CrossDiskDisjoint} applied to the FAL $M_j$, $S$ is disjoint from all standard crossing disks of $M_j$, and so is preserved under the belted sum operations used to reconstruct $M$.  Thus $S$ is a  thrice-punctured sphere in $M$. By Proposition \ref{prop:3psDisjointFromCD} applied to $M_j$, $S$ is orthogonal to the reflection surface in $M_j$ and has two meridional knot circle punctures and one longitudinal crossing circle puncture.   Since the belted sum operation on crossing disks preserves both the reflection surface and component types, $S$ must have these properties in $M$ as well.   The classification of non-standard disks in FALs implies that $S$ is a crossing disk in $M$ as well. But then $M$ contributes two crossing disks to $M_j$, which is a contradiction.  Therefore each summand with at least three crossing circles must be $b$-prime.

Now suppose a summand $M_i$ has two crossing circles.  At least one crossing disk $D$ in $M_i$ is used in the belted sum reconstruction of $M$, which implies that it is the result of slicing along and gluing consecutive crossing disks in $M$.  By definition, the gluing results in a flat crossing circle, so that $D$ must be flat in $M_i$.  Since $M_i$ has two crossing circles, at least one of which is flat, it cannot be the fully twisted Borromean rings, completing the proof.
\end{proof}

It may seem unnatural to ignore singly-separated disks in defining the canonical belted sum decomposition.  Indeed, Lemma \ref{lem:CrossDiskDisjoint} can be generalized to show that all disks involved in any separating pair are disjoint.  One can then slice along all of them and glue corresponding boundary spheres to construct an alternative decomposition.  The components with a singly-separated disk will glue to Whitehead links, and the rest to either $b$-prime FALs with at least three crossing circles or the Borromean rings with at least one flat crossing circle.  Further decomposing Borromean rings summands into two Whitehead link complements in both this decomposition and the canonical one results in the same decomposition.   For simplicity of exposition, we presented the crossing disk method.

Neither approach, that of excluding or of including singly-separated disks, guarantees uniqueness over all belted sum decompositions.  The difficulty is that a non-FAL summand may be produced in the process of decomposing an FAL.  The characterizations of thrice-punctured spheres in Section \ref{sec:3ps} and of separating pairs in Section \ref{sec:SP} cannot be applied to non-FALs, and it is unclear at this writing how to proceed.  For now, then, we content ourselves with the level of uniqueness provided in Theorem \ref{thm:CanonIsBPrime}.

\noindent
\textbf{Acknowledgements}

This research was supported in part by NSF-REU Grants DMS-1461286 and DMS-1758020, as well as California State University, San Bernardino.  We are very grateful to Christian Millichap for many (many!) conversations and suggestions that have significantly improved this work.  Thanks also to Jeff Meyer for helpful conversations.

\end{document}